\documentclass[a4paper]{article}
\usepackage{indentfirst,latexsym,bm,amsfonts,amssymb,amsmath,amsfonts,mathrsfs}
\newtheorem{Definition}{Definition}[section]
\newtheorem{Theorem}{Theorem}[section]
\newtheorem{Example}{Example}[section]
\newtheorem{Lemma}{Lemma}[section]
\newtheorem{Corollary}{Corollary}[section]

\title{\bf Demi-linear Analysis II\linebreak
---Demi-distributions}

\author{Li Ronglu\\ %and Wen Songlong\\%[4pt]
{\it\normalsize Department of Mathematics, Harbin Institute of Technology, Harbin,
150001, China}\\[3pt]
Zhong Shuhui\\
{\it\normalsize Department of Mathematics, Tianjin University, Tianjin,
300072, China}\\{\it\normalsize Corresponding author E-mail: shuhuizhong@126.com}\\[3pt]
Kim Dohan\\ %and Park Siehee\\%[4pt]
{\it\normalsize Department of Mathematics, Seoul National University, Seoul,
151-742, Korea}\\{\it\normalsize Corresponding author E-mail: dhkim@snu.ac.kr
}\\[3pt]
Wu Junde\\ %and Park Siehee\\%[4pt]
{\it\normalsize School of Mathematics Sciences, Zhejiang University, Hangzhou 310027, China}\\{\it\normalsize Corresponding author E-mail: wjd@zju.edu.cn}\\[3pt]
}
\begin{document}
\setcounter{footnote}{1}
\renewcommand{\thefootnote}{\fnsymbol{footnote}}
\footnotetext{
{This work was supported by the National Natural Science
Foundation of China (Grant No. 11126165).}}
\date{}
\maketitle
\def\abstractname{}
\begin{abstract}
\noindent{\bf Abstract.} In this paper, we establish a demi-distributions theory which develops the usual distribution theory, in particular, we show that many conclusions as differentiations, Fourier transforms and convolutions can be generalized to the demi-distributions theory. \vspace{6pt}


\end{abstract}\vspace{20pt}

Let $X$ be a topological vector space and $\mathcal N(X)$ the family
of neighborhoods of $0\in X$, $C(0)$
the set of complex valued functions $\gamma$ satisfying
\begin{enumerate}
\item $\gamma: \mathbb C \rightarrow \mathbb C$;
\item $\lim_{t\rightarrow0}\gamma(t)=\gamma(0)=0$;
\item $|\gamma(t)|\geq|t|$ if $|t|\leq1$.
\end{enumerate}

Let $Y$ be a topological vector space. A mapping $f:X\rightarrow Y$ is said to be {\it demi-linear} if $f(0)=0$
and there exist $\gamma\in C(0)$ and $U\in\mathcal N(X)$ such that
every $x\in X$, $u\in U$ and $t\in\{t\in\mathbb K:|t|\leq1\}$ yield
$r,\ s\in\mathbb K$ for which $|r-1|\leq|\gamma(t)|$,
$|s|\leq|\gamma(t)|$ and $f(x+tu)=rf(x)+sf(u)$.

We denote by ${\mathscr L}_{\gamma,U}(X,Y)$ the demi-linear mappings
related to $\gamma\in C(0)$ and $U\in\mathcal N(X)$, and by
${\mathscr K}_{\gamma,U}(X,Y)$ the subfamily of ${\mathscr L}_{\gamma,U}(X,Y)$ satisfying the following property: if $x\in X,\
u\in U$ and $|t|\leq1$ then
$$f(x+tu)=f(x)+sf(u)$$
for some $s$ with $|s|\leq|\gamma(t)|$.

As was stated in [1, 2], the family of demi-linear mappings is a
natural and valuable extension of the family of linear operators.

For $a>0$, $\mathscr{D}_a=\big\{\xi\in\Bbb C^{\Bbb R^n}:\xi\mbox{ is
infinitely differentiable and }\xi(x)=0\mbox{ whenever
}\linebreak|x|=\sqrt{x_1^2+\cdots+x_n^2}>a\big\}$ has the locally
convex Fr\'echet topology which is given by the norm sequence
$\{\|\xi\|_p=\sup_{|x|\leq a}\max_{|q|\leq
p}|D^q\xi(x)|\}_{p=0}^\infty$. 

Let $\mathscr{D}=\bigcup_{m=1}^\infty\mathscr{D}_m$ be the strict
inductive limit of $\{\mathscr{D}_m\}$.

Let $\mathscr{S}=\{\xi\in\Bbb C^{\Bbb R^n}:\xi\mbox{ is infinitely
differentiable and rapidly decreasing}\}$. With the norm sequence
$\{\|\xi\|_p=\sup_{|k|,|q|\leq p,x\in\mathbb
R^n}|x^kD^q\xi(x)|\}_{p=0}^\infty$, $\mathscr S$ is a locally convex
Fr\'echet space, where $k=(k_1,k_2,\cdots,k_n)$ is a multi-index and
$x^k=x_1^{k_1}x_2^{k_2}\cdots x_n^{k_n}$.

The spaces $\{\mathscr{D}_a,\mathscr{D},\mathscr{S}\}$ are called the space of test functions.

The distributions of L. Schwartz, the generalized distributions of
Beurling and the ultradistributions of Roumieu are continuous linear
functionals defined on some suitable spaces of test functions.

In this paper, $E\in\{\mathscr{D}_a,\mathscr{D},\mathscr{S}\}$ is a
space of test functions and a function $f:E\rightarrow\mathbb{C}$ is
called a demi-distribution if $f$ is continuous and demi-linear. Thus, the family of demi-distributions includes all usual
distributions and many nonlinear functionals [1, 2].

By using the equicontinuity results in [1], we show that the family of demi-linear
mappings can be used to develop the theory of distributions. For
instance, in the case of the usual distributions the simplest
equation $y'=0$ has solutions $y=constant$ only. However, we will
show that the equation $y'=0$ has extremely many solutions which are
nonlinear demi-linear functionals, and the equation $y'=f$ also has
extremely many solutions which are demi-distributions. Moreover, we
will show that the family of demi-distributions is closed with
respect to extremely many of nonlinear transformations such as
$|f(\cdot)|$, $\sin|f(\cdot)|$, $e^{|f(\cdot)|}-1$, etc.

\section{Demi-distributions}

Firstly, note that the space $\mathscr D$ is an $(LF)$ space and so $\mathscr D$ is
both barrelled and bornological [3, Page 222]. Thus, $\mathscr D$ is
$C-sequential$ [3, Page 118]. There is an important fact which says
that every sequentially continuous linear operator from a
$C-sequential$ locally convex space to a locally convex space must
be continuous [3, Page 118]. Now, we we improve the fact as following:

\begin{Theorem}
Let $X,Y$ be locally convex spaces, $U\in\mathcal N(X)$ and
$\gamma_0(t)=t,\ \forall\,t\in\mathbb C$. If $X$ is $C-sequential$
and $f\in\mathscr K_{\gamma_0,U}(X,Y)$ is sequentially continuous,
then $f$ must be continuous.

Proof. Let $V\in\mathcal N(Y)$. Pick balanced convex neighborhoods
$U_0\in\mathcal N(X)$ and $V_0\in\mathcal N(Y)$ such that
$U_0\subset U$, $V_0\subset V$.

Let $W=f^{-1}(V_0)$. For $u,w\in U_0\bigcap W$ and scalars $\alpha,
\beta$ with $|\alpha|+|\beta|\leq1$, it follows from $f\in\mathscr
K_{\gamma_0,U}(X,Y)$ that
$$f(\alpha u+\beta w)=f(\alpha u)+sf(w)=s_1f(u)+sf(w)\in s_1V_0+sV_0,$$
where $|s_1|\leq|\gamma_0(\alpha)|=|\alpha|$,
$|s|\leq|\gamma_0(\beta)|=|\beta|$. Then
$|s_1|+|s|\leq|\alpha|+|\beta|\leq1$ and so $f(\alpha u+\beta w)\in
V_0$, $\alpha u+\beta w\in U_0\bigcap W$. This shows that
$U_0\bigcap W$ is both balanced and convex.

Let $x_k\rightarrow0$ in $X$. Then $x_k\in U_0$ eventually. Since
$f$ is sequentially continuous, $f(x_k)\rightarrow f(0)=0$ and
$f(x_k)\in V_0$ eventually, i.e., $x_k\in W$ eventually. Thus,
$x_k\in U_0\bigcap W$ eventually and so $U_0\bigcap W$ is a
sequential neighborhood of $0\in X$. Since $X$ is $C-sequential$,
$U_0\bigcap W\in\mathcal N(X)$ and $f(U_0\bigcap W)\subset
V_0\subset V$. This shows that $f$ is continuous at $0$.

Suppose $(x_\alpha)_{\alpha\in I}$ is a net in $X$ such that
$x_\alpha\rightarrow x\in X$. Pick an $\alpha_0\in I$ for which
$x_\alpha-x\in U,\ \forall\,\alpha\geq\alpha_0$. For
$\alpha\geq\alpha_0$,
$$f(x_\alpha)=f(x+x_\alpha-x)=f(x)+s_\alpha f(x_\alpha-x),\ |s_\alpha|\leq|\gamma_0(1)|=1.$$
But $f(x_\alpha-x)\rightarrow f(0)=0$ and so $f(x_\alpha)\rightarrow
f(x)$. $\square$
\end{Theorem}

\begin{Theorem}
Suppose that $E\in\{\mathscr D_a,\mathscr D,\mathscr S\}$,
$\gamma\in C(0)$, $U\in\mathcal N(E)$ and $Y$ is a topological
vector space. If $f,f_\nu\in\mathscr L_{\gamma,U}(E,Y)$ are
continuous ($\nu=1,2,3,\cdots$) and $f_\nu(\xi)\rightarrow f(\xi)$
at each $\xi\in E$, then for every bounded $B\subset E$, $\lim_\nu
f_\nu(\xi)=f(\xi)$ uniformly for $\xi\in B$.

Proof. Since both $\mathscr D_a$ and $\mathscr S$ are Fr\'echet
spaces having the Montel property, we only need to consider
$\mathscr D$. Suppose that $B\subset \mathscr D$ is bounded but
$\lim_\nu f_\nu(\xi)=f(\xi)$ is not uniform for $\xi\in B$. Then
there exist $V\in\mathcal N(Y)$, $\{\xi_k\}\subset B$ and integers
$\nu_1<\nu_2<\cdots$ such that
$$f_{\nu_k}(\xi_k)-f(\xi_k)\not\in V,\ k=1,2,3,\cdots.$$

Pick a balanced $W\in\mathcal N(Y)$ for which $W+W+W\subset V$.
Since $B$ is bounded in $\mathscr D$, $B\subset\mathscr D_m$ for
some $m\in\mathbb N$ [3, p.219] and $B$ is relatively compact in the
Fr\'echet space $\mathscr D_m$ [4, Th. 1.6.2]. By passing to a
subsequence if necessary, we say that
$\xi_k\rightarrow\xi\in\mathscr D_m$. Since $f_\nu(\eta)\rightarrow
f(\eta)$ at each $\eta\in\mathscr D_m$, $\{f_\nu\}_1^\infty$ is
pointwise bounded on $\mathscr D_m$ and so $\{f_\nu\}_1^\infty$ is
equicontinuous on $\mathscr D_m$ by Th. 3.1 of [1]. By Cor. 3.1 of
[1], $\lim_k f_\nu(\xi_k)=f_\nu(\xi)$ uniformly for $\nu\in\mathbb
N$ and so there is a $k_0\in\mathbb N$ such that
$f_\nu(\xi_k)-f_\nu(\xi)\in W$ for all $\nu\in\mathbb N$ and
$k>k_0$. Since $f:\mathscr D\rightarrow Y$ is continuous and
$f_\nu(\xi)\rightarrow f(\xi)$, there exist $\nu_0,\ k_1\in\mathbb
N$ such that $f(\xi)-f(\xi_k)\in W$ for all $k>k_1$ and
$f_\nu(\xi)-f(\xi)\in W$ for all $\nu>\nu_0$.

Pick an integer $k_2\geq k_0+k_1$ for which $\nu_k>\nu_0$ whenever
$k>k_2$. Then for every $k>k_2$ we have that
\begin{align*}
f_{\nu_k}(\xi_k)-f(\xi_k)&=f_{\nu_k}(\xi_k)-f_{\nu_k}(\xi)+f_{\nu_k}(\xi)-f(\xi)+f(\xi)-f(\xi_k)\\
&\in W+W+W\subset V.
\end{align*}
This is a contradiction and so $\lim_\nu f_\nu(\xi)=f(\xi)$
uniformly for $\xi\in B$. $\square$
\end{Theorem}

In general, $\mathscr L_{\gamma,U}(X,Y)\subsetneqq\mathscr
W_{\gamma,U}(X,Y)$ where $Y$ is locally convex. Using Th. 4.1 of [1]
instead of Th. 3.1 of [1], the above proof gives an improved result
as follows.

\begin{Theorem}
Let $E\in\{\mathscr D_a,\mathscr D,\mathscr S\}$, $\gamma\in C(0)$
and $U\in\mathcal N(E)$. Let $Y$ be a locally convex space and
$f,f_\nu\in\mathscr W_{\gamma,U}(E,Y)$ are continuous,
$\nu=1,2,3,\cdots$. If $f_\nu(\xi)\rightarrow f(\xi)$ at each
$\xi\in E$, then for every bounded $B\subset E$, $\lim_\nu
f_\nu(\xi)=f(\xi)$ uniformly for $\xi\in B$.
\end{Theorem}

For $\mathscr K_{\gamma_0,U}(X,Y)$, we have a more strong result as
follows.

\begin{Theorem}
Let $E\in\{\mathscr D_a,\mathscr D,\mathscr S\}$, $U\in\mathcal
N(E)$ and $Y$ be a locally convex space. Let $f_\nu\in\mathscr
K_{\gamma_0,U}(E,Y)$ be continuous, $\forall\,\nu\in\mathbb N$. If
$\lim_\nu f_\nu(\xi)=f(\xi)$ exists at each $\xi\in E$, then $f$ is
also a continuous mapping in $\mathscr K_{\gamma_0,U}(E,Y)$ and for
every bounded $B\subset E$, $\lim_\nu f_\nu(\xi)=f(\xi)$ uniformly
for $\xi\in B$.

Proof. Only need to consider $E=\mathscr D$. Let $\xi\in E$,
$\eta\in U$ and $|t|\leq1$. Then $f(\xi+t\eta)=\lim_\nu
f_\nu(\xi+t\eta)=\lim_\nu(f_\nu(\xi)+s_\nu f_\nu(\eta))$, where
$|s_\nu|\leq|\gamma_0(t)|=|t|\leq1$. Say that $s_{\nu_k}\rightarrow
s$. Then $|s|=\lim_k|s_{\nu_k}|\leq|\gamma_0(t)|$ and
$$f(\xi+t\eta)=\lim_kf_{\nu_k}(\xi+t\eta)=\lim_k(f_{\nu_k}(\xi)+s_{\nu_k}f_{\nu_k}(\eta))=f(\xi)+sf(\eta).$$
Thus, $f\in\mathscr K_{\gamma_0,U}(E,Y)$.

Let $\xi_k\rightarrow\xi$ in $\mathscr D$. Then
$\xi_k\rightarrow\xi$ in $\mathscr D_m$ for some $m\in\mathbb N$ [3,
p.219]. Since $f_\nu(\cdot)\rightarrow f(\cdot)$,
$\{f_\nu\}_1^\infty$ is pointwise bounded on $\mathscr D_m$ and, by
Cor. 3.1 of [1], $\lim_kf_\nu(\xi_k)=f_\nu(\xi)$ uniformly for
$\nu\in\mathbb N$. Then $\lim_kf(\xi_k)=\lim_k\lim_\nu
f_\nu(\xi_k)=\lim_\nu\lim_kf_\nu(\xi_k)=\lim_\nu f_\nu(\xi)=f(\xi)$.
Thus, $f:\mathscr D\rightarrow Y$ is sequentially continuous and so
$f$ is continuous by Th. 1.1.

Now the desired follows from Th. 1.2. $\square$
\end{Theorem}

Henceforth, $E\in\{\mathscr D_a,\mathscr D,\mathscr S\}$.

\begin{Definition}
$f:E\rightarrow\mathbb C$ is called a demi-distribution if $f$ is
continuous and $f\in\mathscr L_{\gamma,U}(E,\mathbb C)$ for some
$\gamma\in C(0)$ and $U\in\mathcal N(E)$.
\end{Definition}

Let $E^{(\gamma,U)}$ be the family of demi-distributions related to
$\gamma\in C(0)$ and $U\in\mathcal N(E)$. Let $[E^{(\gamma,U)}]$ be
the span $(E^{(\gamma,U)})$ in $\mathbb C^E$, i.e.,
$[E^{(\gamma,U)}]=\big\{\mbox{finite sum }\sum t_kf_k:t_k\in\mathbb
C,f_k\in E^{(\gamma,U)}\big\}$.

Let $E'$ be the space of usual distributions, i.e., $E'$ is the
space of continuous linear functionals. Obviously, $E'\subset
E^{(\gamma,U)},\ \forall\,U\in\mathcal N(E),\ \gamma\in C(0)$.

\begin{Example}
(1) For every $f\in L_{loc}^1(\Bbb R^n)$ define $[f]:\mathscr
D\rightarrow\Bbb R$ by $[f](\xi)=\int_{\mathbb
R^n}|f(x)\xi(x)|\,dx,\ \xi\in\mathscr D$.

Let $\gamma\in C(0)$ and $\xi,\eta\in \mathscr D$, $|t|\leq1$. For
every $x\in\Bbb R^n$ there exists $\alpha(x)\in[-|t|,|t|]$ such that
$|\xi(x)+t\eta(x)|=|\xi(x)|+\alpha(x)|\eta(x)|$ and
\begin{align*}
[f](\xi+t\eta)&=\int_{\mathbb R^n}|f(x)(\xi+t\eta)(x)|\,dx\\
&=\int_{\mathbb R^n}|f(x)||\xi(x)+t\eta(x)|\,dx\\
&=\int_{\mathbb R^n}|f(x)|[|\xi(x)|+\alpha(x)|\eta(x)|]\,dx\\
&=\int_{\mathbb R^n}|f(x)\xi(x)|\,dx+\int_{\mathbb
R^n}\alpha(x)|f(x)\eta(x)|\,dx\\
&=[f](\xi)+\int_{\mathbb R^n}\alpha(x)|f(x)\eta(x)|\,dx.
\end{align*}

If $\int_{\mathbb R^n}|f(x)\eta(x)|\,dx=0$, then
$0\leq|\int_{\mathbb
R^n}\alpha(x)|f(x)\eta(x)|\,dx|\leq\int_{\mathbb
R^n}|\alpha(x)f(x)\eta(x)|dx\linebreak\leq|t|\int_{\mathbb
R^n}|f(x)\eta(x)|\,dx=0$ and so $ \int_{\mathbb
R^n}\alpha(x)|f(x)\eta(x)|\,dx=0=0[f](\eta)$, where
$0\leq|\gamma(t)|$.

If $\int_{\mathbb R^n}|f(x)\eta(x)|\,dx\neq0$, then
$|\frac{\int_{\mathbb R^n}\alpha(x)|f(x)\eta(x)|\,dx}{\int_{\mathbb
R^n}|f(x)\eta(x)|\,dx}|\leq|t|\leq|\gamma(t)|$ and so \linebreak
$\int_{\mathbb R^n}\alpha(x)|f(x)\eta(x)|\,dx=s\int_{\mathbb
R^n}|f(x)\eta(x)|\,dx=s[f](\eta)$ where $s=\frac{\int_{\mathbb
R^n}\alpha(x)|f(x)\eta(x)|\,dx}{\int_{\mathbb
R^n}|f(x)\eta(x)|\,dx}$, $|s|\leq|t|\leq|\gamma(t)|$. Thus,
$$[f](\xi+t\eta)=[f](\xi)+s[f](\eta),\ |s|\leq|\gamma(t)|,$$
i.e., $[f]\in\mathscr K_{\gamma,\mathscr D}(\mathscr D,\mathbb
R)\bigcap\mathscr D^{(\gamma,\mathscr D)}$ but $[f]$ is not a usual
distribution.

(2) Let $\mathscr D_1(\mathbb R)=\big\{\xi\in\mathbb R^\mathbb
R:\xi\mbox{ is infinitely differentiable and }\xi(x)=0\mbox{ for
}|x|>1\big\}$. Let $\gamma(t)=\frac{\pi}{2}t$ for $t\in\mathbb R$
and $U=\big\{\xi\in\mathscr D_1(\mathbb R):\max_{|x|\leq
1}|\xi(x)|<1\big\}$. Define $f:\mathscr D_1(\mathbb
R)\rightarrow\mathbb R$ by
$$f(\xi)=\int_{-\infty}^\infty|\sin \xi(x)|\,dx,\ \xi\in\mathscr D_1(\mathbb R).$$

It is easy to show that if $a\in\mathbb R$ and $u,t\in[-1,1]$ then
$\sin (a+tu)=\sin a+s\sin u$ with $|s|\leq\frac{\pi}{2}|t|$. Hence,
for $\xi\in\mathscr D_1(\mathbb R),\ \eta\in U$ and $|t|\leq1$ we
have that
\begin{align*}
f(\xi+t\eta)&=\int_{-\infty}^\infty|\sin[\xi(x)+t\eta(x)]|\,dx\\
&=\int_{-\infty}^\infty|\sin\xi(x)+\alpha(x)\sin\eta(x)|\,dx\qquad (|\alpha(x)|\leq\frac{\pi}{2}|t|)\\
&=\int_{-\infty}^\infty[|\sin\xi(x)|+\beta(x)|\sin\eta(x)|]\,dx\qquad (|\beta(x)|\leq|\alpha(x)|\leq\frac{\pi}{2}|t|)\\
&=\int_{-\infty}^\infty|\sin\xi(x)|\,dx+\int_{-\infty}^\infty\beta(x)|\sin\eta(x)|\,dx\\
&=\int_{-\infty}^\infty|\sin\xi(x)|\,dx+s\int_{-\infty}^\infty|\sin\eta(x)|\,dx\qquad(|s|\leq\frac{\pi}{2}|t|=|\gamma(t)|)\\
&=f(\xi)+sf(\eta),\quad|s|\leq\frac{\pi}{2}|t|=|\gamma(t)|.
\end{align*}
Thus, $f\in\mathscr K_{\gamma,U}(\mathscr D_1(\mathbb R),\mathbb
R)\bigcap(\mathscr D_1(\mathbb R))^{(\gamma,U)}$ but $f$ is not a
usual distribution.

(3) For the case of $\mathbb R^n=\mathbb R$, we write that $\mathscr
S=\mathscr S(\mathbb R)$. Let $U=\big\{\eta\in\mathscr S(\mathbb
R):\sup_{x\in\mathbb R}|\eta(x)|<1\big\}$ and $\gamma(t)=et$ for
$t\in\mathbb C$. Then define $g:\mathscr S(\mathbb
R)\rightarrow\mathbb C$ by
$$g(\xi)=\sqrt{-1}\int_{-1}^1(e^{|\xi(x)|}-1)\,dx,\ \xi\in\mathscr S(\mathbb R).$$

For $\xi\in\mathscr S(\mathbb R),\eta\in U$ and $|t|\leq1$,
\begin{align*}
g(\xi+t\eta)&=\sqrt{-1}\int_{-1}^1(e^{|\xi(x)+t\eta(x)|}-1)\,dx\\
&=\sqrt{-1}\int_{-1}^1(e^{|\xi(x)|+\alpha(x)|\eta(x)|}-1)\,dx\qquad(\alpha(x)\in[-|t|,|t|])\\
&=\sqrt{-1}\int_{-1}^1(e^{|\xi(x)|+\alpha(x)|\eta(x)|}-e^{\alpha(x)|\eta(x)|}+e^{\alpha(x)|\eta(x)|}-1)\,dx\\
&=\sqrt{-1}\int_{-1}^1
e^{\alpha(x)|\eta(x)|}(e^{|\xi(x)|}-1)\,dx+\sqrt{-1}\int_{-1}^1(e^{\alpha(x)|\eta(x)|}-1)\,dx.
\end{align*}
If $\int_{-1}^1(e^{|\xi(x)|}-1)\,dx=0$, then $0\leq\int_{-1}^1
e^{\alpha(x)|\eta(x)|}(e^{|\xi(x)|}-1)\,dx\leq
e^{|t|}\int_{-1}^1(e^{|\xi(x)|}-1)\,dx=0$ and so
$\sqrt{-1}\int_{-1}^1
e^{\alpha(x)|\eta(x)|}(e^{|\xi(x)|}-1)\,dx=0=\sqrt{-1}\int_{-1}^1(e^{|\xi(x)|}-1)\,dx=g(\xi)=rg(\xi)$
where $r=1$, $|r-1|=0\leq|\gamma(t)|$. If $\int_{-1}^1
(e^{|\xi(x)|}-1)\,dx\neq0$, then $\sqrt{-1}\int_{-1}^1
e^{\alpha(x)|\eta(x)|}(e^{|\xi(x)|}-1)\,dx=\frac{\int_{-1}^1
e^{\alpha(x)|\eta(x)|}(e^{|\xi(x)|}-1)\,dx}{\int_{-1}^1
(e^{|\xi(x)|}-1)\,dx}\sqrt{-1}\int_{-1}^1(e^{|\xi(x)|}-1)\,dx$,
where
\begin{align*}
&|\frac{\int_{-1}^1
e^{\alpha(x)|\eta(x)|}(e^{|\xi(x)|}-1)\,dx}{\int_{-1}^1(e^{|\xi(x)|}-1)\,dx}-1|
=\frac{|\int_{-1}^1
(e^{\alpha(x)|\eta(x)|}-1)(e^{|\xi(x)|}-1)\,dx|}{\int_{-1}^1(e^{|\xi(x)|}-1)\,dx}\\
=&\frac{|\int_{-1}^1
e^{\theta(x)\alpha(x)|\eta(x)|}\alpha(x)|\eta(x)|(e^{|\xi(x)|}-1)\,dx|}{\int_{-1}^1(e^{|\xi(x)|}-1)\,dx}\qquad(0\leq\theta(x)\leq1)\\
\leq&\frac{\int_{-1}^1
e^{|t|}|t|(e^{|\xi(x)|}-1)\,dx}{\int_{-1}^1(e^{|\xi(x)|}-1)\,dx}\qquad(\because|\alpha(x)|\leq|t|,\
|\eta(x)|\leq1)\\
=&e^{|t|}|t|\leq e|t|=|\gamma(t)|.\qquad(\because|t|\leq1)
\end{align*}
Thus, $\sqrt{-1}\int_{-1}^1
e^{\alpha(x)|\eta(x)|}(e^{|\xi(x)|}-1)\,dx=r\sqrt{-1}\int_{-1}^1(e^{|\xi(x)|}-1)\,dx=rg(\xi)$
where $|r-1|\leq|\gamma(t)|$.

If $g(\eta)=\sqrt{-1}\int_{-1}^1(e^{|\eta(x)|}-1)\,dx=0$, then
$\eta(x)=0$ a.e. in $[-1,1]$ and
$g(\xi+t\eta)=\sqrt{-1}\int_{-1}^1(e^{|\xi(x)+t\eta(x)|}-1)\,dx=\sqrt{-1}\int_{-1}^1(e^{|\xi(x)|}-1)\,dx
=g(\xi)=rg(\xi)+sg(\eta)$ where $r=1$ and $s=0,\
|r-1|=0\leq|\gamma(t)|$, $|s|=0\leq|\gamma(t)|$.

Suppose that
$g(\eta)=\sqrt{-1}\int_{-1}^1(e^{|\eta(x)|}-1)\,dx\neq0$. Then
$$\sqrt{-1}\int_{-1}^1(e^{\alpha(x)|\eta(x)|}-1)\,dx
=\frac{\int_{-1}^1(e^{\alpha(x)|\eta(x)|}-1)\,dx}{\int_{-1}^1(e^{|\eta(x)|}-1)\,dx}g(\eta),$$
where
\begin{align*}
|\frac{\int_{-1}^1(e^{\alpha(x)|\eta(x)|}-1)\,dx}{\int_{-1}^1(e^{|\eta(x)|}-1)\,dx}|
&=\frac{|\int_{-1}^1
e^{\delta(x)\alpha(x)|\eta(x)|}\alpha(x)|\eta(x)|\,dx|}{\int_{-1}^1
e^{\theta(x)|\eta(x)|}|\eta(x)|\,dx}\quad(0\leq\delta(x),\
\theta(x)\leq1)\\
&\leq\frac{|\int_{-1}^1
e^{\delta(x)\alpha(x)|\eta(x)|}|\alpha(x)||\eta(x)|\,dx|}{\int_{-1}^1|\eta(x)|\,dx}\\
&\leq \frac{\int_{-1}^1
e^{|t|}|t||\eta(x)|\,dx}{\int_{-1}^1|\eta(x)|\,dx}\\
&=e^{|t|}|t|\leq e|t|=|\gamma(t)|.%\qquad(\because|\eta(x)|\leq1,\
%|\alpha(x)|\leq|t|)
\end{align*}
Then $g(\xi+t\eta)=rg(\xi)+sg(\eta)$ where $|r-1|\leq|\gamma(t)|$,
$|s|\leq|\gamma(t)|$, i.e., $g\in\mathscr L_{\gamma,U}(\mathscr
S(\mathbb R),\mathbb C)$. Since $\xi_k\rightarrow\xi$ in $\mathscr
S$ implies that $\|\xi_k-\xi\|_0=\sup_{x\in\mathbb
R}|\xi_k(x)-\xi(x)|\rightarrow 0$ and so
$g(\xi_k)=\sqrt{-1}\int_{-1}^1(e^{|\xi_k(x)|}-1)\,dx\rightarrow\sqrt{-1}\int_{-1}^1(e^{|\xi(x)|}-1)\,dx=g(\xi)$,
i.e., $g:\mathscr S(\mathbb R)\rightarrow\mathbb C$ is continuous.
Thus, $g\in (\mathscr S(\mathbb R))^{(\gamma,U)}$.
\end{Example}

For every $C\geq1$ and $\varepsilon>0$, $\mathscr
K_{C,\varepsilon}(\mathbb R,\mathbb R)$ includes a lot of nonlinear
functions. Pick a $h\in\mathscr K_{C,\varepsilon}(\mathbb R,\mathbb
R)$ and let $f(x+iy)=ih(|x+iy|),\ \forall\,x+iy\in\mathbb C$. If
$|u+iv|<\varepsilon$ and $|t|\leq1$, then
$f[x+iy+t(u+iv)]=ih(|x+iy+t(u+iv)|)=ih(|x+iy|+\alpha|u+iv|)=ih(|x+iy|)+sih(|u+iv|)=f(x+iy)+sf(u+iv)$,
where $\alpha\in[-|t|,|t|]\subset[-1,1]$ and $|s|\leq C|\alpha|\leq
C|t|$. This shows that $f\in\mathscr K_{C,\varepsilon}(\mathbb
C,\mathbb C)$ and, therefore, $\mathscr K_{C,\varepsilon}(\mathbb
C,\mathbb C)$ also includes a lot of nonlinear functions.

Let $E^{[\gamma,U]}=\{f\in\mathscr K_{\gamma,U}(E,\mathbb C):f\mbox{
is continuous}\}$. Then $E^{[\gamma,U]}\subset E^{(\gamma,U)}$.

\begin{Theorem}
If $A\subset E'$ is an equicontinuous family of distributions and
$\varepsilon>0$, then there is a $U\in\mathcal N(E)$ such that
$$\big\{h\circ f:h\in\mathscr L_{\gamma,\varepsilon}(\mathbb C,\mathbb C),f\in A\}\subset E^{(\gamma,U)},\
\forall\,\gamma\in C(0),$$
$$\big\{h\circ f:h\in\mathscr K_{\gamma,\varepsilon}(\mathbb C,\mathbb C),f\in A\}
\subset E^{[\gamma,U]},\ \forall\,\gamma\in C(0).$$

Proof. Since $A$ is equicontinuous, there is a $U\in\mathcal N(E)$
such that $|f(\eta)|<\varepsilon,\ \forall\,f\in A,\eta\in U$. Let
$\xi\in E,\ \eta\in U$ and $|t|\leq1$. For $h\in\mathscr
L_{\gamma,\varepsilon}(\mathbb C,\mathbb C)$ and $f\in A$,
\begin{align*}
(h\circ f)(\xi+t\eta)&=h(f(\xi)+tf(\eta))\\
&=r(h\circ f)(\xi)+s(h\circ f)(\eta),\ \ \ |r-1|\leq|\gamma(t)|,\
|s|\leq|\gamma(t)|.
\end{align*}
Thus, $h\circ f\in\mathscr L_{\gamma,U}(E,\mathbb C)$.

Suppose that $h\in\mathscr L_{\gamma,\varepsilon}(\mathbb C,\mathbb
C)$ and $w_k\rightarrow w$ in $\mathbb C$. Then
\begin{align*}
\lim_kh(w_k)&=\lim_kh(w+w_k-w)=\lim_kh(w+\frac{2(w_k-w)}{\varepsilon}\frac{\varepsilon}{2})\\
&=\lim_k[r_kh(w)+s_kh(\frac{\varepsilon}{2})],
\end{align*}
where
$|r_k-1|\leq|\gamma(\frac{2(w_k-w)}{\varepsilon})|\rightarrow0$ and
$|s_k|\leq|\gamma(\frac{2(w_k-w)}{\varepsilon})|\rightarrow0$, i.e.,
$r_k\rightarrow1,\ s_k\rightarrow0$. Thus, $h(w_k)\rightarrow h(w)$,
$h$ is continuous. But $A\subset E'$ and so $h\circ
f:E\rightarrow\mathbb C$ is continuous for $h\in\mathscr
L_{\gamma,\varepsilon}(\mathbb C,\mathbb C)$ and $f\in A$. $\square$
\end{Theorem}

\begin{Corollary}
If $f\in E'$ is a usual distribution and $\varepsilon>0$, then there
is a $U\in\mathcal N(E)$ such that
$$\big\{h\circ f:h\in\mathscr L_{\gamma,\varepsilon}(\mathbb C,\mathbb C)\}\subset E^{(\gamma,U)},\
\forall\,\gamma\in C(0),$$
$$\big\{h\circ f:h\in\mathscr K_{\gamma,\varepsilon}(\mathbb C,\mathbb C)\}\subset E^{[\gamma,U]},\
\forall\,\gamma\in C(0).$$
\end{Corollary}

\begin{Corollary}
If $A\subset E'$ is a pointwise bounded family of usual
distributions, then for every $\varepsilon>0$ there is a
$U\in\mathcal N(E)$ such that
$$\big\{h\circ f:h\in\mathscr L_{\gamma,\varepsilon}(\mathbb C,\mathbb C),f\in A\}\subset E^{(\gamma,U)},\
\forall\,\gamma\in C(0),$$
$$\big\{h\circ f:h\in\mathscr K_{\gamma,\varepsilon}(\mathbb C,\mathbb C),f\in A\}\subset E^{[\gamma,U]},\
\forall\,\gamma\in C(0).$$

Proof. Both $\mathscr D_a$ and $\mathscr S$ are Fr\'echet spaces. So
for the case of $E\in\{\mathscr D_a,\mathscr S\}$, $A$ is
equicontinuous by Th. 3.1 of [1].

Since $\mathscr D$ is an $(LF)$ space, $\mathscr D$ is barrelled [3,
p.222]. Then $A$ is also equicontinuous for the case of $E=\mathscr
D$ [3, Th. 9.3.4]. $\square$
\end{Corollary}

\begin{Corollary}
For every $U\in\mathcal N(E)$ and the polar $U^\circ=\{f\in
E':|f(\eta)|\leq1,\ \forall\,\eta\in U\}$,
$$\{h\circ f:h\in\mathscr L_{\gamma,1}(\mathbb C,\mathbb C),f\in U^\circ\}\subset E^{(\gamma,U)},\ \forall\,\gamma\in C(0),$$
$$\{h\circ f:h\in\mathscr K_{\gamma,1}(\mathbb C,\mathbb C),f\in U^\circ\}\subset E^{[\gamma,U]},\ \forall\,\gamma\in C(0).$$

Proof. $U^\circ$ is equicontinuous [3, p.129]. $\square$
\end{Corollary}

\begin{Theorem}
Let $\gamma_1\in C(0)$ for which $\sup_{|t|\leq1}|\gamma_1(t)|=1$,
$|\gamma_1(\alpha)|\leq|\gamma_1(\beta)|$ whenever
$|\alpha|\leq|\beta|\leq1$, e.g., $\gamma_1(t)=\sqrt{|t|}$. For
every $U\in\mathcal N(E)$, $f\in E^{[\gamma_1,U]}$ and
$\varepsilon>0$ there is a $V\in\mathcal N(E)$ such that
$\gamma_1\circ\gamma_1\in C(0)$ and
$$h\circ f\in E^{(\gamma_1\circ\gamma_1, V)},\ \forall\,h\in\mathscr L_{\gamma_1,\varepsilon}(\mathbb C,\mathbb C),$$
$$h\circ f\in E^{[\gamma_1\circ\gamma_1, V]},\ \forall\,h\in\mathscr K_{\gamma_1,\varepsilon}(\mathbb C,\mathbb C).$$

Proof. Pick a $W\in\mathcal N(E)$ for which $|f(\eta)|<\varepsilon,\
\forall\,\eta\in W$. Let $h\in\mathscr
L_{\gamma_1,\varepsilon}(\mathbb C,\mathbb C)$, $\xi\in E$, $\eta\in
V=U\bigcap W$ and $|t|\leq1$. Then
\begin{align*}
(h\circ f)(\xi+t\eta)&=h(f(\xi+t\eta))=h(f(\xi)+\alpha
f(\eta))\qquad
(|\alpha|\leq|\gamma_1(t)|\leq1)\\
&=rh(f(\xi))+sh(f(\eta))=r(h\circ f)(\xi)+s(h\circ f)(\eta),
\end{align*}
where $|r-1|\leq|\gamma_1(\alpha)|\leq|\gamma_1(\gamma_1(t))|$,
$|s|\leq|\gamma_1(\alpha)|\leq|\gamma_1(\gamma_1(t))|$.

As in the proof of Th. 1.5, $h$ is continuous and so $h\circ f\in
E^{(\gamma_1\circ\gamma_1,V)}$.

Similarly, $h\circ f\in E^{[\gamma_1\circ\gamma_1,V]}$ whenever
$h\in\mathscr K_{\gamma_1,\varepsilon}(\mathbb C,\mathbb C)$.
$\square$
\end{Theorem}

\begin{Example}
(1) Let $h(z)=|z|,\ \forall\,z\in\mathbb C$. Then $h\in\mathscr
K_{\gamma_0,\mathbb C}(\mathbb C,\mathbb C)$ where $\gamma_0(t)=t$.
Let $U\in\mathcal N(E)$ and $\gamma_1\in C(0)$ as in Th. 1.6. Then
for every $f\in E^{(\gamma_1,U)}$ and $a>0$ there is a
$V_a\in\mathcal N(E)$ such that $V_a\subset U$ and $|f(\eta)|<a,\
\forall\,\eta\in V_a$.

Let $a>0$, $\xi\in E$, $\eta\in V_a$ and $|t|\leq1$. Then
$$|f(\xi+t\eta)|=|rf(\xi)+\alpha f(\eta)|=|r||f(\xi)|+s|f(\eta)|,$$
where $||r|-1|\leq|r-1|\leq|\gamma_1(t)|$,
$|s|\leq|\alpha|\leq|\gamma_1(t)|$. Thus,
$\gamma_0\circ\gamma_1=\gamma_1$ and $|f(\cdot)|=h\circ f\in
E^{(\gamma_1,V_a)},\ \forall\,a>0$.

(2) Let $\gamma_1(t)=\sqrt{|t|}$, $\gamma_2(t)=\frac{\pi}{2}\,t$,
$\forall\,t\in\mathbb C$. Let $U\in\mathcal N(E)$ and $f\in
E^{[\gamma_1,U]}$. There is a $V\in\mathcal N(E)$ such that
$V\subset U$ and $|f(\eta)|<1$, $\forall\,\eta\in V$. Define
$\sin|f(\cdot)|:E\rightarrow\mathbb C$ by
$\sin|f(\cdot)|(\xi)=\sin|f(\xi)|,\ \xi\in E$. For $\xi\in E$,
$\eta\in V$ and $|t|\leq1$,
\begin{align*}
\sin|f(\cdot)|(\xi+t\eta)&=\sin|f(\xi+t\eta)|=\sin|f(\xi)+\alpha
f(\eta)|\quad\ (|\alpha|\leq|\gamma_1(t)|=\sqrt{|t|}\leq1)\\
&=\sin[|f(\xi)|+\beta|f(\eta)|]\quad\ (|\beta|\leq|\alpha|\leq1)\\
&=\sin|f(\xi)|+s\sin|f(\eta)|\
(|s|\leq\frac{\pi}{2}|\beta|\leq\frac{\pi}{2}|\alpha|\leq\frac{\pi}{2}\sqrt{|t|}=|(\gamma_2\circ\gamma_1)(t)|)\\
&=\sin|f(\cdot)|(\xi)+s\sin|f(\cdot)|(\eta).
\end{align*}
Thus, $\gamma_2\circ\gamma_1\in C(0)$ and $\sin|f(\cdot)|\in
E^{[\gamma_2\circ\gamma_1,V]}$.

(3) If $h(z)=e^{|z|}-1,\ \forall\,z\in\mathbb C$,
$\gamma_1(t)=\sqrt{|t|}$ and $\gamma(t)=e^2t$, then $h\in\mathscr
L_{\gamma,1}(\mathbb C,\mathbb C)$ and for every $f\in
E^{[\gamma_1,U]}$ there is a $V\in\mathcal N(E)$ such that
$e^{|f(\cdot)|}-1=h\circ f\in E^{(\gamma\circ\gamma_1,V)}$.
\end{Example}

Even $f$ is a nonzero usual distribution, each of $|f(\cdot)|$,
$\sin|f(\cdot)|$ and $e^{|f(\cdot)|}-1$ can not be a usual
distribution. However, Th. 1.6 shows that the family of
demi-distributions is closed with respect to infinitely many of
nonlinear transformations.

Henceforth, in the notations $E^{(\gamma,U)}$ and $E^{[\gamma,U]}$
we always confess that $\gamma\in C(0)$ and $U\in\mathcal N(E)$.

\begin{Definition}
$f_k\stackrel{w\ast}{\longrightarrow}f$ in $E^{(\gamma,U)}$ means
that $f_k,f\in E^{(\gamma,U)}$ for all $k\in\mathbb N$ and
$f_k(\xi)\rightarrow f(\xi)$ at each $\xi\in E$, and $f_k\rightarrow
f$ in $E^{(\gamma,U)}$ means that $f_k,f\in E^{(\gamma,U)}$ for all
$k\in\mathbb N$ and for every bounded $B\subset E$,
$\lim_kf_k(\xi)=f(\xi)$ uniformly for $\xi\in B$.
\end{Definition}

Now Th. 1.2 gives the following

\begin{Theorem}
$f_k\rightarrow f$ in $E^{(\gamma,U)}$ if and only if
$f_k\stackrel{w\ast}{\longrightarrow}f$ in $E^{(\gamma,U)}$.
\end{Theorem}

\begin{Definition}
A sequence $\{f_k\}\subset E^{(\gamma,U)}$ (resp., $E^{[\gamma,U]}$)
is $w\ast$ Cauchy if $\lim_kf_k(\xi)$ exists at each $\xi\in E$.
$E^{(\gamma,U)}$ (resp., $E^{[\gamma,U]}$) is said to be
sequentially complete if for every $w\ast$ Cauchy sequence $\{f_k\}$
in $E^{(\gamma,U)}$ (resp., $E^{[\gamma,U]}$) there exists $f\in
E^{(\gamma,U)}$ (resp., $E^{[\gamma,U]}$) such that $f_k\rightarrow
f$, i.e., for every bounded $B\subset E$, $\lim_kf_k(\xi)=f(\xi)$
uniformly for $\xi\in B$.
\end{Definition}

\begin{Theorem}
Both $\mathscr D_a^{(\gamma,U)}$ and $\mathscr S^{(\gamma, V)}$ are
sequentially complete for every $\gamma\in C(0)$, $U\in\mathcal
N(\mathscr D_a)$ and $V\in\mathcal N(\mathscr S)$. Moreover,
$\mathscr D^{[\gamma_0,W]}$ is also sequentially complete for
$\gamma_0(t)=t$ and $W\in\mathcal N(\mathscr D)$.

Proof. Let $E\in\{\mathscr D_a,\mathscr S\}$, $U\in\mathcal N(E)$
and $\gamma\in C(0)$. If $\{f_k\}\subset E^{(\gamma,U)}$ and
$\lim_kf_k(\xi)=f(\xi)$ exists at each $\xi\in E$, then $\{f_k\}$ is
equicontinuous by Th. 3.1 of [1]. If $\xi_\nu\rightarrow\xi$ in $E$,
then $\lim_\nu f_k(\xi_\nu)=f_k(\xi)$ uniformly for $k\in\mathbb N$
and $\lim_\nu f(\xi_\nu)=\lim_\nu\lim_kf_k(\xi_\nu)=\lim_k\lim_\nu
f_k(\xi_\nu)=\lim_kf_k(\xi)=f(\xi)$. Thus, $f:E\rightarrow\mathbb C$
is continuous.

Let $\xi\in E,\ \eta\in U$ and $|t|\leq1$. Then
$f(\xi+t\eta)=\lim_kf_k(\xi+t\eta)=\lim_k[r_kf_k(\xi)+s_kf_k(\eta)]$
where $|r_k-1|\leq|\gamma(t)|$, $|s_k|\leq|\gamma(t)|$. By passing
to a subsequence if necessary, we assume that $r_k\rightarrow r$ and
$s_k\rightarrow s$. Then $|r-1|\leq|\gamma(t)|$,
$|s|\leq|\gamma(t)|$ and
$f(\xi+t\eta)=\lim_k[r_kf_k(\xi)+s_kf_k(\eta)]=rf(\xi)+sf(\eta)$,
$f\in\mathscr L_{\gamma,U}(E,\mathbb C)$. Thus, $f\in
E^{(\gamma,U)}$.

Now $f_k\stackrel{w\ast}{\longrightarrow}f$ in $E^{(\gamma,U)}$ and
so $f_k\rightarrow f$ in $E^{(\gamma,U)}$ by Th. 1.7.

The completeness of $\mathscr D^{[\gamma_0,W]}$ follows from Th.
1.4. $\square$
\end{Theorem}

For $E\in\{\mathscr D,\mathscr S\}$ and $G\subset\mathbb R^n$, let
$E_G=\{\xi\in E:supp\,\xi\subset G\}$. Each $f\in E^{(\gamma,U)}$
yields $f|_{E_G}:E_G\rightarrow\mathbb C$ by $f|_{E_G}(\xi)=f(\xi),\
\forall\,\xi\in E_G$. Then $E_\emptyset=\{0\}$ and
$f|_{E_{\emptyset}}=0,\ \forall\,f\in E^{(\gamma,U)}$.

\begin{Definition}
$E=\mathscr D$. For $f\in E^{(\gamma,U)}$ let
$$supp\,f=\mathbb R^n\backslash\,\Big[\,\bigcup\{G\subset\mathbb R^n:G\mbox{ is open, }f|_{E_G}=0\}\Big].$$
\end{Definition}

\begin{Theorem}
$E=\mathscr D$. For every $f\in E^{(\gamma,U)}$ there is an open
$G_0\subset\mathbb R^n$ such that $f|_{E_{G_0}}=0$ and
$supp\,f=\mathbb R^n\backslash G_0$.

Proof. Let $G_0=\mathbb R^n\backslash supp\,f$ and
$\{G_\alpha:\alpha\in I\}=\{G\subset\mathbb R^n:G\mbox{ is open,
}f|_{E_G}=0\}$. Then $G_0=\bigcup_{\alpha\in I}G_{\alpha}$ is open
and $supp\,f=\mathbb R^n\backslash G_0$.

Suppose that $f|_{E_{G_0}}\neq0$. There is a $\xi\in E$ such that
$supp\,\xi\subset G_0$ but $f(\xi)\neq0$. Then $\mathbb
R^n\backslash supp\,\xi\supset\mathbb R^n\backslash G_0=supp\,f$ and
$\mathbb R^n=(\mathbb R^n\backslash
supp\,\xi)\bigcup(\bigcup_{\alpha\in I}G_\alpha)$. By the partition
of unity, there is a sequence $\{\xi_k\}\subset\mathscr D$ such that
$\sum_{k=1}^\infty\xi_k(x)=1$ for all $x\in\mathbb R^n$ and each
$supp\,\xi_k\subset(\mathbb R^n\backslash supp\,\xi)$ or some
$G_\alpha$, and each $x\in supp\,\xi$ has a neighborhood which
intersects finitely many of $supp\,\xi_k$ only. But $supp\,\xi$ is
compact and so there is an open $G\subset\mathbb R^n$ such that
$supp\,\xi\subset G$ and $G$ intersects finitely many of
$supp\,\xi_k$ only. Hence, $\xi\xi_k=0$ for all but finitely many of
$k'$s. Say that $\{k:\xi\xi_k\neq0\}=\{1,2,\cdots,m\}$. Then
$\xi(x)=\sum_{k=1}^m\xi(x)\xi_k(x),\ \forall\,x\in\mathbb R^n$. For
$k\leq m$, $\xi\xi_k\neq0$ shows that $supp\,\xi_k\not\subset\mathbb
R^n\backslash supp\,\xi$ and so $supp\xi_k\subset G_\alpha$ for some
$\alpha\in I$. Thus, $f(\xi\xi_k)=0,\ k=1,2,\cdots,m$.

Pick a $p\in\mathbb N$ such that $\frac{1}{p}\xi\xi_k\in U,\
k=1,2,\cdots,m$. Since $supp\,(\frac{1}{p}\xi\xi_k)\subset
supp\,\xi_k\subset G_\alpha$ for some $\alpha\in I$,
$f(\frac{1}{p}\xi\xi_k)=0,\ k=1,2,\cdots,m$. Then
\begin{align*}
f(\xi)&=f(\sum_{k=1}^m\xi\xi_k)=f(\sum_{k=1}^{m-1}\xi\xi_k+(p-1)\frac{1}{p}\xi\xi_m+\frac{1}{p}\xi\xi_m)\\
&=r_1f(\sum_{k=1}^{m-1}\xi\xi_k+(p-1)\frac{1}{p}\xi\xi_m)+s_1f(\frac{1}{p}\xi\xi_m)
\quad(|r_1-1|\leq|\gamma(1)|,\ |s_1|\leq|\gamma(1)|)\\
&=r_1f(\sum_{k=1}^{m-1}\xi\xi_k+(p-1)\frac{1}{p}\xi\xi_m)\\
&\qquad\cdots\quad\cdots\\
&=r_1r_2\cdots r_pf(\sum_{k=1}^{m-1}\xi\xi_k)\\
&=r_1r_2\cdots
r_pf(\sum_{k=1}^{m-2}\xi\xi_k+(p-1)\frac{1}{p}\xi\xi_{m-1}+\frac{1}{p}\xi\xi_{m-1})\\
&=r_1\cdots r_pr_{p+1}\cdots r_{2p}f(\sum_{k=1}^{m-2}\xi\xi_k)\\
&\qquad\cdots\quad\cdots\\
&=(\prod_{\nu=1}^{mp-1}r_\nu)f(\frac{1}{p}\xi\xi_1)=0.
\end{align*}

This contradicts that $f(\xi)\neq0$. Hence, $f|_{E_{G_0}}=0$.
$\square$
\end{Theorem}

\begin{Corollary}
For $E=\mathscr D$ and $f\in E^{(\gamma,U)}$, $f=0$ if and only if
$supp\,f=\emptyset$.

Proof. If $f=0$, then $f|_{E_{\mathbb R^n}}=0$ and
$supp\,f\subset\mathbb R^n\backslash\mathbb R^n=\emptyset$, i.e.,
$supp\,f=\emptyset$.

By Th. 1.9, $supp\,f=\mathbb R^n\backslash G_0$ for some open
$G_0\subset\mathbb R^n$ and $f|_{E_{G_0}}=0$. If
$supp\,f=\emptyset$, then $\emptyset=\mathbb R^n\backslash G_0$ and
so $G_0=\mathbb R^n$ and $f|_{E_{\mathbb R^n}}=0$, i.e.,
$f=f|_{E_{\mathbb R^n}}=0$. $\square$
\end{Corollary}

\begin{Corollary}
Let $f,g\in E^{(\gamma,U)}$. Then $f=g$ if and only if for every
$x\in\mathbb R^n$ there is a neighborhood $G$ of $x$ such that
$f|_{E_G}=g|_{E_G}$.

Proof. If for every $x\in\mathbb R^n$ there is an open
$G_x\subset\mathbb R^n$ such that $x\in G_x$ and
$f|_{E_{G_x}}=g|_{E_{G_x}}$, i.e., $(f-g)|_{E_{G_x}}=0$, then
$supp\,(f-g)\subset\mathbb R^n\backslash(\bigcup_{x\in\mathbb
R^n}G_x)=\emptyset$ and so $f-g=0$. $\square$
\end{Corollary}

\section{Differentiation}
$E\in\{\mathscr D_a,\mathscr D,\mathscr S\}$,\ \
$[E^{(\gamma,U)}]=span\,(E^{(\gamma,U)})$ in $\mathbb C^E$.

\begin{Definition}
Let $f=\sum_{k=1}^m\alpha_kf_k\in[E^{(\gamma,U)}]$ where each
$\alpha_k\in\mathbb C$, $f_k\in E^{(\gamma,U)}$. Observing each
$\xi\in E$ is a function defined on $\mathbb R^n$, for
$j\in\{1,2,\cdots,n\}$ define $\frac{\partial f}{\partial
x_j}:E\rightarrow\mathbb C$ by
$$\frac{\partial f}{\partial x_j}(\xi)=f(-\frac{\partial\xi}{\partial x_j}),\ \xi\in E.$$
Then $\frac{\partial}{\partial
x_j}(\sum_{k=1}^m\alpha_kf_k)(\xi)=(\sum_{k=1}^m\alpha_kf_k)(-\frac{\partial\xi}{\partial
x_j})=\sum_{k=1}^m\alpha_kf_k(-\frac{\partial\xi}{\partial
x_j})=\linebreak\sum_{k=1}^m\alpha_k\frac{\partial f_k}{\partial
x_j}(\xi)=(\sum_{k=1}^m\alpha_k\frac{\partial f_k}{\partial
x_j})(\xi)$ and so $\frac{\partial}{\partial
x_j}(\sum_{k=1}^m\alpha_kf_k)=\sum_{k=1}^m\alpha_k\frac{\partial
f_k}{\partial x_j}$ for
$\sum_{k=1}^m\alpha_kf_k\in[E^{(\gamma,U)}]$.
\end{Definition}

For a multi-index $\alpha=(\alpha_1,\cdots,\alpha_n)$ let
$|\alpha|=\alpha_1+\cdots+\alpha_n$ and
$D^\alpha=\frac{\partial^{|\alpha|}}{\partial
x_1^{\alpha_1}\cdots\partial x_n^{\alpha_n}}$. As in the case of
usual distributions, for $f\in[E^{(\gamma,U)}]$ and every
multi-index $\alpha=(\alpha_1,\cdots,\alpha_n)$, $(D^\alpha
f)(\xi)=f((-1)^{|\alpha|}D^\alpha\xi),\ \forall\,\xi\in E$.
Evidently, we have that $\frac{\partial^2f}{\partial x_i\partial
x_j}=\frac{\partial^2f}{\partial x_j\partial x_i}$,
$\frac{\partial^5f}{\partial x_1\partial x_2^2\partial
x_3^2}=\frac{\partial^5f}{\partial x_3^2\partial x_1\partial
x_2^2}$, etc.

\begin{Lemma}
For every multi-index $\alpha$, $D^\alpha:E\rightarrow E$ is a
continuous linear operator.

Proof. For $E=\mathscr D_a$ or $\mathscr S$, the conclusion is
obvious.

Let $\xi_k\rightarrow0$ in $\mathscr D$. Then
$\{\xi_k\}\subset\mathscr D_m$ for some $m\in\mathbb N$ and
$\xi_k\rightarrow 0$ in $\mathscr D_m$ since $\mathscr D_m$ is a
subspace of $\mathscr D$ [3, p.219]. Then
$\|D^\alpha\xi_k\|_p\leq\|\xi_k\|_{|\alpha|+p}$ for all $p\in\mathbb
N$ and so $D^\alpha\xi_k\rightarrow0$ in $\mathscr D_m$, i.e.,
$D^\alpha\xi_k\rightarrow0$ in $\mathscr D$. Thus,
$D^\alpha:\mathscr D\rightarrow\mathscr D$ is sequentially
continuous. Then $D^\alpha:\mathscr D\rightarrow\mathscr D$ is
continuous since $\mathscr D$ is bornological and $C-sequential$.
See also Th. 1.1. $\square$
\end{Lemma}

\begin{Theorem}
Let $\alpha$ be a multi-index. For every $U\in\mathcal N(E)$ there
is a $V\in\mathcal N(E)$ such that
$$\{D^\alpha f:f\in
E^{(\gamma,U)}\}\subset E^{(\gamma,V)},\ \forall\,\gamma\in C(0),$$
$$\{D^\alpha f:f\in
[E^{(\gamma,U)}]\}\subset[E^{(\gamma,V)}],\ \forall\,\gamma\in
C(0).$$

Moreover, if $f_k, f\in E^{(\gamma,U)}$ and
$f_k\stackrel{w\ast}{\longrightarrow}f$, i.e., $f_k(\xi)\rightarrow
f(\xi)$ at each $\xi\in E$, then for every bounded $B\subset E$,
$\lim_k(D^\alpha f_k)(\xi)=(D^\alpha f)(\xi)$ uniformly for $\xi\in
B$.

Proof. Let $U\in\mathcal N(E)$. By Lemma 2.1, there is a
$V\in\mathcal N(E)$ for which $(-1)^{|\alpha|}D^\alpha\eta\in U$,
$\forall\,\eta\in V$.

Let $f\in E^{(\gamma,U)}$, $\xi\in E$, $\eta\in V$ and $|t|\leq1$.
Then
\begin{align*}
(D^\alpha f)(\xi+t\eta)&=f((-1)^{|\alpha|}D^\alpha\xi+t(-1)^{|\alpha|}D^\alpha\eta)\\
&=rf((-1)^{|\alpha|}D^\alpha\xi)+sf((-1)^{|\alpha|}D^\alpha\eta)\\
&=r(D^\alpha f)(\xi)+s(D^\alpha f)(\eta),
\end{align*}
where $|r-1|\leq|\gamma(t)|,\ |s|\leq|\gamma(t)|$. Thus, $D^\alpha
f\in\mathscr L_{\gamma,V}(E,\mathbb C)$.

Since both $(-1)^{|\alpha|}D^\alpha:E\rightarrow E$ and
$f:E\rightarrow\mathbb C$ are continuous, $D^\alpha
f=f\circ(-1)^{|\alpha|}D^\alpha:E\rightarrow \mathbb C$ is also
continuous and so $D^\alpha f\in E^{(\gamma,V)}$.

Suppose that $f_k,f\in E^{(\gamma,U)},\ f_k(\xi)\rightarrow f(\xi),\
\forall\,\xi\in E$, and $B\subset E$ is bounded. Then
$(-1)^{|\alpha|}D^\alpha(B)=\{(-1)^{|\alpha|}D^\alpha\xi:\xi\in B\}$
is bounded and, by Th. 1.7, $\lim_k(D^\alpha
f_k)(\xi)=\lim_kf_k((-1)^{|\alpha|}D^\alpha\xi)=f((-1)^{|\alpha|}D^\alpha\xi)=(D^\alpha
f)(\xi)$ uniformly for $\xi\in B$. $\square$
\end{Theorem}

\begin{Corollary}
Let $\alpha$ be a multi-index and $U\in\mathcal N(E)$. There is a
$V\in\mathcal N(E)$ such that for every $\gamma\in C(0)$ the
differentiation operator
$D^\alpha:[E^{(\gamma,U)}]\rightarrow[E^{(\gamma,V)}]$ is linear and
$w\ast-w\ast$ continuous.

Proof. For $f,g\in[E^{(\gamma,U)}]$ and $t\in\mathbb C$,
$(D^\alpha(f+tg))(\xi)=(f+tg)((-1)^{|\alpha|}D^\alpha\xi)=f((-1)^{|\alpha|}D^\alpha\xi)+tg((-1)^{|\alpha|}D^\alpha\xi)
=(D^\alpha f)(\xi)+t(D^\alpha g)(\xi),\ \forall\,\xi\in E$, i.e.,
$D^\alpha(f+tg)=D^\alpha f+tD^\alpha g$.

If $f_k,f\in[E^{(\gamma,U)}]$ such that
$f_k\stackrel{w\ast}{\longrightarrow}f$, then
$$(D^\alpha f_k)(\xi)=f_k((-1)^{|\alpha|}D^\alpha\xi)\rightarrow f((-1)^{|\alpha|}D^\alpha\xi)=(D^\alpha f)(\xi),
\ \forall\,\xi\in E.\ \square$$
\end{Corollary}

\begin{Example}
(1) Let $f\in L_{loc}^1(\mathbb R^n)$, $\gamma\in C(0)$ and
$$[f](\xi)=\int_{\mathbb R^n}|f(x)\xi(x)|\,dx,\ \xi\in\mathscr D.$$
Then $[f]\in\mathscr D^{[\gamma,\mathscr D]}$ (see Exam. 1.1(1)),
and
$$(D^\alpha [f])(\xi)=\int_{\mathbb R^n}|f(x)(D^\alpha\xi)(x)|\,dx,\ \forall\,\xi\in\mathscr D.$$

(2) $\gamma$ and $U$ as in Exam. 1.1(2), and
$$f(\xi)=\int_{-\infty}^\infty|\sin\xi(x)|\,dx,\ \forall\,\xi\in\mathscr D_1(\mathbb R).$$
Then $f\in(\mathscr D_1(\mathbb R))^{[\gamma,U]}$ and
\begin{align*}
(D^\alpha f)(\xi)&=\int_{\infty}^\infty|\sin[(-1)^{|\alpha|}(D^\alpha\xi)(x)]|\,dx\\
&=\int_{-\infty}^\infty|\sin(D^\alpha\xi)(x)|\,dx,\
\forall\,\xi\in\mathscr D_1(\mathbb R).
\end{align*}

(3) Let $U\in\mathcal N(E)$ and $\gamma(t)=\sqrt{|t|},\
\forall\,t\in\mathbb C$. For every $f\in E^{[\gamma,U]}$, both
$\sin|f(\cdot)|$ and $e^{|f(\cdot)|}-1$ are demi-distributions, see
Exam. 1.2. Then for every multi-index $\alpha$,
$$D^\alpha\sin|f(\cdot)|=\sin|D^\alpha f(\cdot)|,\ \ D^\alpha(e^{|f(\cdot)|}-1)=e^{|D^\alpha f(\cdot)|}-1.$$
In general, Th. 1.6 shows that $h\circ f$ is a demi-distribution for
every $h\in\mathscr L_{\gamma,\varepsilon}(\mathbb C,\mathbb C)$.
Then
$$D^\alpha(h\circ f)=h\circ D^\alpha f.$$
In fact, $(D^\alpha(h\circ f))(\xi)=(h\circ
f)((-1)^{|\alpha|}D^\alpha\xi)=h[f((-1)^{|\alpha|}D^\alpha\xi)]=h[(D^\alpha
f)(\xi)]=(h\circ D^\alpha f)(\xi),\ \forall\,\xi\in E$.

(4) $E\in\{\mathscr D_a,\mathscr S\}$ and
$\{\|\cdot\|_p\}_{p=0}^\infty$ is the usual norm sequence on $E$.
For $p\in\mathbb N$ and $\varepsilon>0$, let
$U_{p,\varepsilon}=\{\eta\in E:\|\eta\|_p<\varepsilon\}$. Then for
every multi-index $\alpha$ and $\varepsilon>0$,
$$D^\alpha f\in E^{(\gamma,U_{p+|\alpha|,\varepsilon})},
\ \forall\,f\in E^{(\gamma,U_{p,\varepsilon})},\ \gamma\in C(0),\
p\in\mathbb N,$$
$$D^\alpha f\in [E^{(\gamma,U_{p+|\alpha|,\varepsilon})}],
\ \forall\,f\in [E^{(\gamma,U_{p,\varepsilon})}],\ \gamma\in C(0),\
p\in\mathbb N.$$ In fact, $\|\eta\|_{p+|\alpha|}<\varepsilon$
implies
$\|(-1)^{|\alpha|}D^\alpha\eta\|_p\leq\|\eta\|_{p+|\alpha|}<\varepsilon$.
\end{Example}

\begin{Definition}
$\zeta:\mathbb R^n\rightarrow\mathbb C$ is called a multiplier in
$E$ if for every $\xi\in E$ the pointwise product $\zeta\xi\in E$
and $\zeta\xi_k\rightarrow0$ in $E$ whenever $\xi_k\rightarrow0$.
For a multiplier $\zeta$ in $E$ and $f\in[E^{(\gamma,U)}]$, define
$\zeta f:E\rightarrow\mathbb C$ by $(\zeta f)(\xi)=f(\zeta\xi),\
\forall\,\xi\in E$.
\end{Definition}

\begin{Theorem}
If $\zeta$ is a multiplier in $E$, then for every $U\in\mathcal
N(E)$ there is a $V\in\mathcal N(E)$ such that
$$\{\zeta f:f\in E^{(\gamma,U)}\}\subset E^{(\gamma,V)},\ \forall\,\gamma\in C(0),$$
$$\{\zeta f:f\in[E^{(\gamma,U)}]\}\subset[E^{(\gamma,V)}],\ \forall\,\gamma\in C(0).$$

Proof. The correspondence $\xi\mapsto\zeta\xi$ is a continuous
linear operator from $E$ into $E$ and so there is a $V\in\mathcal
N(E)$ such that $\zeta\eta\in U$ for all $\eta\in V$.

Let $f\in E^{(\gamma,U)}$ and $\xi\in E$, $\eta\in V$, $|t|\leq1$.
Then
\begin{align*}
(\zeta f)(\xi+t\eta) &=f(\zeta\xi+t\zeta\eta)=rf(\zeta\xi)+sf(\zeta\eta)\\
&=r(\zeta f)(\xi)+s(\zeta f)(\eta),\ \ |r-1|\leq|\gamma(t)|,\
|s|\leq|\gamma(t)|.
\end{align*}
Thus, $\zeta f\in \mathscr L_{\gamma,V}(E,\mathbb C)$. The
continuity of $\zeta f:E\rightarrow\mathbb C$ follows from the
continuity of $f$ and the continuity of the correspondence
$\xi\mapsto\zeta\xi$. Hence, $\zeta f\in E^{(\gamma,V)}$. $\square$
\end{Theorem}

\begin{Lemma}
Let $E\in\{\mathscr D_a(\mathbb R),\mathscr D(\mathbb R),\mathscr
S(\mathbb R)\}$ be a space of test functions defined on $\mathbb R$,
i.e., $n=1$. Pick a $\zeta\in E$ for which
$\int_{-\infty}^{\infty}\zeta(x)\,dx=1$ and define $A:E\rightarrow
E$ by $A(\xi)=\xi-(\int_{-\infty}^\infty\xi(x)\,dx)\zeta,\ \xi\in
E$. Then $A$ is a continuous linear operator, $\int_{-\infty}^\infty
A(\xi)(x)\,dx=0$ for all $\xi\in E$ and $A(\xi^{(k)})=\xi^{(k)},\
\forall\,\xi\in E,\ k\in\mathbb N$.

Proof. For $\xi,\eta\in E$ and $t\in\mathbb C$,
$A(\xi+t\eta)=\xi+t\eta-(\int_{-\infty}^\infty(\xi+t\eta)(x)\,dx)\zeta
=\xi-(\int_{-\infty}^\infty\xi(x)\,dx)\zeta+t\eta-t(\int_{-\infty}^\infty\eta(x)\,dx)\zeta=A(\xi)+tA(\eta).$

Since $1\in E'$, if $\xi_k\rightarrow\xi$ in $E$, then
$\int_{-\infty}^\infty\xi_k(x)\,dx\rightarrow\int_{-\infty}^\infty\xi(x)\,dx$,
$A(\xi_k)=\xi_k-(\int_{-\infty}^\infty\xi_k(x)\,dx)\zeta
\rightarrow\xi-(\int_{-\infty}^\infty\xi(x)\,dx)\zeta=A(\xi)$ and so
$A$ is sequentially continuous. Since $E$ is bornological, $A$ is
continuous.

For $\xi\in E$ and $k\geq1$,
$A(\xi^{(k)})=\xi^{(k)}-(\int_{-\infty}^\infty\xi^{(k)}(x)\,dx)\zeta=\xi^{(k)}$.
$\square$
\end{Lemma}

For usual distributions, the equation $y'=0$ has solutions $y=const$
only. However, for demi-distributions in $E^{(\gamma,U)}$, the
equation $y'=0$ has extremely many solutions which are not
constants, and the equation $y'=f$ also has extremely many solutions
which are demi-distributions.

\begin{Lemma}
Let $E$ be a space of test functions defined on $\mathbb R$. Let
$\gamma\in C(0)$ and $U\in\mathcal N(E)$. For $y\in E^{(\gamma,U)}$,
$y'=0$ if and only if $y(\xi)=0$ whenever
$\int_{-\infty}^\infty\xi(x)\,dx=0$.

Proof. Suppose that $y'=0$ and $\xi\in E$ for which
$\int_{-\infty}^\infty\xi(x)\,dx=0$. Letting
$\eta(x)=\int_{-\infty}^x\xi(t)\,dt$, $\eta\in E$ and $\xi=\eta'$.
Then $y(\xi)=y(-(-\eta)')=y'(-\eta)=0$.

The converse is obvious. $\square$
\end{Lemma}

\begin{Example}
(1) Let $\gamma\in C(0)$ and $U\in\mathcal N(\mathscr D(\mathbb
R))$. Pick a nonzero $\xi_0\in\mathscr D(\mathbb R)$, e.g.,
$\xi_0(x)=e^{1/(x^2-1)}$ if $|x|<1$ and $0$ if $|x|\geq1$. Since the
functional $\xi\mapsto\int_{-\infty}^\infty\xi(x)\,dx$ is continuous
on $\mathscr D(\mathbb R)$, i.e., $1\in(\mathscr D(\mathbb R))'$,
there is a $V\in\mathcal N(\mathscr D(\mathbb R))$ such that
$(\int_{-\infty}^\infty\eta(x)\,dx)\xi_0\in U,\ \forall\,\eta\in V$.
Then pick a nonzero $f_0\in(\mathscr D(\mathbb R))^{(\gamma,U)}$ and
define $f:\mathscr D(\mathbb R)\rightarrow\mathbb C$ by
$f(\xi)=f_0((\int_{-\infty}^\infty\xi(x)\,dx)\xi_0),\ \xi\in\mathscr
D(\mathbb R)$.

Evidently, $f$ is continuous. For $\xi\in\mathscr D(\mathbb R),\
\eta\in V$ and $|t|\leq1$,
\begin{align*}
f(\xi+t\eta)&=f_0[(\int_{-\infty}^\infty\xi(x)\,dx)\xi_0+t(\int_{-\infty}^\infty\eta(x)\,dx)\xi_0]\\
&=rf_0((\int_{-\infty}^\infty\xi(x)\,dx)\xi_0)+sf_0((\int_{-\infty}^\infty\eta(x)\,dx)\xi_0)\\
&=rf(\xi)+sf(\eta),\ \ |r-1|\leq|\gamma(t)|,\ |s|\leq|\gamma(t)|.
\end{align*}
Thus, $f\in(\mathscr D(\mathbb R))^{(\gamma,V)}$.

If $\xi\in E$ and $\int_{-\infty}^\infty\xi(x)\,dx=0$, then
$f(\xi)=f_0(0)=0$ so $f'=0$ by Lemma 2.3, i.e., $f$ is a solution of
the equation $y'=0$. If
$\xi_0(x)=\begin{cases}e^{1/(x^2-1)},&|x|<1,\\0,&|x|\geq1\end{cases}$
and $f_0=[1]$, i.e., $f_0(\xi)=\int_{-\infty}^\infty|\xi(x)|\,dx,\
\forall\,\xi\in\mathscr D(\mathbb R)$, then
\begin{align*}
f(\xi)&=f_0[(\int_{-\infty}^\infty\xi(x)\,dx)\xi_0]
=\int_{-\infty}^\infty|(\int_{-\infty}^\infty\xi(t)\,dt)\xi_0(x)|\,dx\\
&=|\int_{-\infty}^\infty\xi(x)\,dx|\int_{-1}^1e^{1/(x^2-1)}\,dx,\
\forall\,\xi\in\mathscr D(\mathbb R).
\end{align*}
Note that this $f$ is not a usual distribution because
\begin{align*}
f(\xi+\eta)&=|\int_{-\infty}^\infty\xi(x)\,dx+\int_{-\infty}^\infty\eta(x)\,dx|\int_{-\infty}^\infty e^{1/(x_2-1)}\,dx\\
&\neq[|\int_{-\infty}^\infty\xi(x)\,dx|+
|\int_{-\infty}^\infty\eta(x)\,dx|]\int_{-\infty}^\infty e^{1/(x^2-1)}\,dx\\
&=f(\xi)+f(\eta)
\end{align*}
for many pairs $\xi,\eta\in\mathscr D(\mathbb R)$, e.g., if
$\xi=\xi_0$, $\eta=-\xi_0$, then $f(\xi+\eta)=f(0)=0$ but
$f(\xi)+f(\eta)=2[\int_{-1}^1e^{1/(x^2-1)}\,dx]^2>0$. But $f'=0$.

The solution
$f(\xi)=f_0[(\int_{-\infty}^\infty\xi(x)\,dx)\xi_0]=|\int_{-\infty}^\infty\xi(x)\,dx|\int_{-\infty}^\infty|\xi_0(x)|\,dx$
determined by $\xi_0\in E$ and $f_0=[1]$ has a nice splitting
property: $f\in(\mathscr D(\mathbb R))^{[\gamma_0,\mathscr D(\mathbb
R)]}$. In fact, for $\xi,\eta\in\mathscr D(\mathbb R)$ and
$t\in\mathbb R$,
\begin{align*}
f(\xi+t\eta)&=|\int_{-\infty}^\infty[\xi(x)+t\eta(x)]\,dx|\int_{-\infty}^\infty|\xi_0(x)|\,dx\\
&=|\int_{-\infty}^\infty\xi(x)\,dx+t\int_{-\infty}^\infty\eta(x)\,dx|\int_{-\infty}^\infty|\xi_0(x)|\,dx\\
&=[|\int_{-\infty}^\infty\xi(x)\,dx|+s|\int_{-\infty}^\infty\eta(x)\,dx|]\int_{-\infty}^\infty|\xi_0(x)|\,dx\\
&=f(\xi)+sf(\eta),\ \ |s|\leq|t|.
\end{align*}

(2) Let $\gamma(t)=\frac{\pi}{2}\,t$ for $t\in\mathbb R$,
$U=\{\xi\in\mathscr D_1(\mathbb R):\max_{|x|\leq1}|\xi(x)|<1\}$ and
$$f_0(\xi)=\int_{-\infty}^\infty|\sin\,\xi(x)|\,dx,\ \xi\in\mathscr D_1(\mathbb R).$$
Then $f_0\in(\mathscr D_1(\mathbb R))^{[\gamma,U]}$ (Exam. 1.1(2)).

Pick a nonzero $\xi_0\in\mathscr D_1(\mathbb R)$ and let
$f(\xi)=f_0[(\int_{-\infty}^\infty\xi(x)\,dx)\xi_0]
=\linebreak\int_{-\infty}^\infty|\sin\,[(\int_{-\infty}^\infty\xi(t)\,dt)\xi_0(x)]|\,dx$
for $\xi\in\mathscr D_1(\mathbb R)$. Then pick a $V\in\mathcal
N(\mathscr D_1(\mathbb R))$ such that
$(\int_{-\infty}^\infty\eta(x)\,dx)\xi_0\in U,\ \forall\,\eta\in V$.

Let $\xi\in\mathscr D_1(\mathbb R)$, $\eta\in V$ and $|t|\leq1$.
Then
$\max_{|x|\leq1}|(\int_{-\infty}^\infty\eta(s)\,ds)\xi_0(x)|\leq1$,
and
\begin{align*}
f(\xi+t\eta)&=\int_{-\infty}^\infty|\sin\,[(\int_{-\infty}^\infty\xi(r)\,dr)\xi_0(x)
+t(\int_{-\infty}^\infty\eta(r)\,dr)\xi_0(x)]|\,dx\\
&=\int_{-\infty}^\infty|\sin\,[(\int_{-\infty}^\infty\xi(r)\,dr)\xi_0(x)]
+\alpha(x)\sin\,[(\int_{-\infty}^\infty\eta(r)\,dr)\xi_0(x)]|\,dx\\
&=\int_{-\infty}^\infty|\sin\,[(\int_{-\infty}^\infty\xi(r)\,dr)\xi_0(x)]|\,dx
+\int_{-\infty}^\infty\beta(x)|\sin\,[(\int_{-\infty}^\infty\eta(r)\,dr)\xi_0(x)]|\,dx,
\end{align*}
where $|\beta(x)|\leq|\alpha(x)|\leq\frac{\pi}{2}|t|$ for all
$x\in\mathbb R$. It is similar to Exam. 1.1,
\begin{align*}
f(\xi+t\eta)&=\int_{-\infty}^\infty|\sin\,[(\int_{-\infty}^\infty\xi(r)\,dr)\xi_0(x)]|\,dx
+s\int_{-\infty}^\infty|\sin\,[(\int_{-\infty}^\infty\eta(r)\,dr)\xi_0(x)]|\,dx\\
&=f(\xi)+sf(\eta),\ \ |s|\leq\frac{\pi}{2}|t|=|\gamma(t)|.
\end{align*}
Then $f\in(\mathscr D_1(\mathbb R))^{[\gamma,V]}$ and $f'=0$.
\end{Example}

In general, we have

\begin{Theorem}
Let $E\in\{\mathscr D_a(\mathbb R),\mathscr D(\mathbb R),\mathscr
S(\mathbb R)\}$, a space of test functions defined on $\mathbb R$.
Let $U\in\mathcal N(E)$, $\gamma\in C(0)$. Then for every $\xi_0\in
E$ and $f_0\in E^{(\gamma,U)}$ there is a $V\in\mathcal N(E)$ such
that the equation $y'=0$ has a solution $f\in E^{(\gamma,V)}$ which
is given by $f(\xi)=f_0[(\int_{-\infty}^\infty\xi(x)\,dx)\xi_0]$ for
$\xi\in E$. If $f_0=1$, then $f_0\in E'$ and
$f(\xi)=f_0[(\int_{-\infty}^\infty\xi(x)\,dx)\xi_0]=(\int_{-\infty}^\infty\xi(x)\,dx)f_0(\xi_0)
=(\int_{-\infty}^\infty\xi(x)\,dx)(\int_{-\infty}^\infty\xi_0(\tau)\,d\tau)
=\linebreak\int_{-\infty}^\infty(\int_{-\infty}^\infty\xi_0(\tau)\,d\tau)\xi(x)\,dx$,
$\forall\,\xi\in E$, i.e.,
$f=\int_{-\infty}^\infty\xi_0(\tau)\,d\tau\in E'$, a constant which
is a usual solution of the equation $y'=0$.
\end{Theorem}

The solutions of the equation $y'=0$ have an interesting property as
follows.

\begin{Theorem}
Let $E\in\{\mathscr D_a(\mathbb R),\mathscr D(\mathbb R)\}$, a space
of test functions defined on $\mathbb R$. Let $U\in\mathcal N(E)$,
$\gamma\in C(0)$ and $y\in E^{[\gamma,U]}$. If $y'=0$, then for
every $\zeta\in E$ with $\int_{-\infty}^\infty\zeta(x)\,dx=1$,
$$y(\xi)=y[(\int_{-\infty}^\infty\xi(x)\,dx)\zeta],\ \forall\,\xi\in E,$$
i.e., if $\zeta_1,\zeta_2\in E$ such that
$\int_{-\infty}^\infty\zeta_1(x)\,dx=\int_{-\infty}^\infty\zeta_2(x)\,dx=1$,
then
$y[(\int_{-\infty}^\infty\xi(x)\,dx)\zeta_1]\\=y[(\int_{-\infty}^\infty\xi(x)\,dx)\zeta_2]=y(\xi)$
for all $\xi\in E$ and, in particular,
$$y(\xi)=y(\eta)\mbox{ whenever }\int_{-\infty}^\infty\xi(x)\,dx=\int_{-\infty}^\infty\eta(x)\,dx=1,$$
$$y(\frac{\xi}{\int_{-\infty}^\infty\xi(x)\,dx})=y(\frac{\eta}{\int_{-\infty}^\infty\eta(x)\,dx})
\mbox{ whenever }\int_{-\infty}^\infty\xi(x)\,dx\neq0\mbox{ and
}\int_{-\infty}^\infty\eta(x)\,dx\neq0.$$

Proof. If $\zeta\in E$ such that $\zeta\neq\xi',\ \forall\,\xi\in
E$, i.e., $\int_{-\infty}^\infty\zeta(x)\neq0$, then
$\frac{1}{\int_{-\infty}^\infty\zeta(x)\,dx}\zeta\in E$ and
$\int_{-\infty}^\infty\frac{\zeta(x)}{\int_{-\infty}^\infty\zeta(t)\,dt}\,dx=1$.
Pick a $\zeta\in E$ for which $\int_{-\infty}^\infty\zeta(x)\,dx=1$
and let $A(\xi)=\xi-(\int_{-\infty}^\infty\xi(x)\,dx)\zeta$ for
$\xi\in E$. By Lemma 2.2, $A:E\rightarrow E$ is a continuous linear
operator and $\int_{-\infty}^\infty A(\xi)(x)\,dx=0,\
\forall\,\xi\in E$. Moreover,
$$A(\xi)(x)=\Big(\int_{-\infty}^x A(\xi)(t)\,dt\Big)',\ \forall\,\xi\in E,\ x\in\mathbb R.$$

Let $\xi\in E$ and pick a $p\in\mathbb N$ for which
$\frac{1}{p}A(\xi)\in U$. Then
\begin{align*}
y(\xi)&=y[(\int_{-\infty}^\infty\xi(x)\,dx)\zeta+A(\xi)]\\
&=y[(\int_{-\infty}^\infty\xi(x)\,dx)\zeta+(p-1)\frac{1}{p}A(\xi)+\frac{1}{p}A(\xi)]\\
&=y[(\int_{-\infty}^\infty\xi(x)\,dx)\zeta+(p-1)\frac{1}{p}A(\xi)]+s_1y(\frac{1}{p}A(\xi))\\
&\qquad\qquad\qquad\cdots\\
&=y[(\int_{-\infty}^\infty\xi(x)\,dx)\zeta]+sy(\frac{1}{p}A(\xi)).
\end{align*}
But
$\frac{1}{p}A(\xi)(x)=(\frac{1}{p}\int_{-\infty}^xA(\xi)(t)\,dt)'$
for all $x\in\mathbb R$ and so $y(\frac{1}{p}A(\xi))=\linebreak
y[-(-\frac{1}{p}\int_{-\infty}^xA(\xi)(t)\,dt)']=y'(-\frac{1}{p}\int_{-\infty}^xA(\xi)(t)\,dt)=0$
since $y'=0$. Therefore,
$$y(\xi)=y[(\int_{-\infty}^\infty\xi(x)\,dx)\zeta],\ \forall\,\xi\in E.$$

If
$\int_{-\infty}^\infty\xi(x)\,dx=\int_{-\infty}^\infty\eta(x)\,dx=1$,
then
$y(\xi)=y[(\int_{-\infty}^\infty\xi(x)\,dx)\zeta]=y(\zeta)=y[(\int_{-\infty}^\infty\eta(x)\,dx)\zeta]=y(\eta)$.
$\square$
\end{Theorem}

For $E\in\{\mathscr D_a(\mathbb R),\mathscr D(\mathbb R)\}$ let
$$E_1=\{\xi\in E:\int_{-\infty}^\infty\xi(x)\,dx=1\}.$$
If $y\in E'$ is a usual distribution such that $y'=0$, then $y$ must
be a constant $C\in\mathbb R$, i.e., $y(\xi)=\int_{-\infty}^\infty
C\xi(x)\,dx$ for all $\xi\in E$. Hence, $y(\xi)=C,\ \forall\,\xi\in
E_1$. Th. 2.3 shows that the same fact holds for the case of
$E^{[\gamma,U]}$.

\begin{Corollary}
$E\in\{\mathscr D_a(\mathbb R),\mathscr D(\mathbb R)\}$,
$U\in\mathcal N(E)$ and $\gamma\in C(0)$. If $y\in E^{[\gamma,U]}$
such that $y'=0$, then $y(\cdot)$ is an invariant on $E_1$, i.e.,
there is a $C\in\mathbb R$ such that
$$y(\xi)=C,\ \forall\,\xi\in E_1.$$
\end{Corollary}

Although Th. 2.3 gives a lot of various solutions of the equation
$y'=0$ for the case of $E^{(\gamma,U)}$ but Th. 2.3 did not give all
solutions. However, for the case of $E^{[\gamma,U]}$ we can give all
solutions of $y'=0$.

\begin{Theorem}
Let $E\in\{\mathscr D_a(\mathbb R),\mathscr D(\mathbb R)\}$ be a
space of test functions defined on $\mathbb R$, $U\in\mathcal N(E)$
and $\gamma\in C(0)$. Then for every $\xi_0\in E$ and $f_0\in
E^{[\gamma,U]}$ there is a $V\in\mathcal N(E)$ such that the
equation $y'=0$ has a solution $f\in E^{[\gamma,V]}$ which is given
by $f(\xi)=f_0[(\int_{-\infty}^\infty\xi(x)\,dx)\xi_0]$ for $\xi\in
E$. Conversely, if $f\in E^{[\gamma,U]}$ is a solution of the
equation $y'=0$, then there exist $\xi_0\in E$ and $f_0\in
E^{[\gamma,U]}$ such that
$f(\xi)=f_0[(\int_{-\infty}^\infty\xi(x)\,dx)\xi_0],\
\forall\,\xi\in E$.

Proof. Let $\xi_0\in E$, $f_0\in E^{[\gamma,U]}$. There is a
$V\in\mathcal N(E)$ such that
$(\int_{-\infty}^\infty\eta(x)\,dx)\xi_0\linebreak\in U$ for all
$\eta\in V$. Let
$f(\xi)=f_0[(\int_{-\infty}^\infty\xi(x)\,dx)\xi_0]$ for $\xi\in E$.
If $\xi\in E$, $\eta\in V$ and $|t|\leq1$, then
$f(\xi+t\eta)=f_0[(\int_{-\infty}^\infty(\xi+t\eta)(x)\,dx)\xi_0]
=f_0[(\int_{-\infty}^\infty\xi(x)\,dx)\xi_0+t(\int_{-\infty}^\infty\eta(x)\,dx)\xi_0]
=f_0[(\int_{-\infty}^\infty\xi(x)\,dx)\xi_0]+sf_0[(\int_{-\infty}^\infty\eta(x)\,dx)\xi_0]=f(\xi)+sf(\eta)$,
where $|s|\leq|\gamma(t)|$. Thus, $f\in E^{[\gamma,V]}$ and $f'=0$:
$$f'(\xi)=f(-\xi')=f_0[(\int_{-\infty}^\infty-\xi'(x)\,dx)\xi_0]=f_0(0)=0,\ \forall\,\xi\in E.$$

Conversely, suppose that $f\in E^{[\gamma,U]}$ and $f'=0$. Pick a
$\,\xi_0\in E$ with $\int_{-\infty}^\infty\xi_0(x)\,dx\linebreak=1$,
and let $f_0=f$. By Th. 2.4,
$$f(\xi)=f[(\int_{-\infty}^\infty\xi(x)\,dx)\xi_0]
=f_0[(\int_{-\infty}^\infty\xi(x)\,dx)\xi_0],\ \forall\,\xi\in E.\
\square$$
\end{Theorem}

We now consider the equation $y'=f$ where $f\in E^{(\gamma,U)}$.

\begin{Theorem}
Let $E\in\{\mathscr D_a(\mathbb R),\mathscr D(\mathbb R)\}$ be a
space of test functions defined on $\mathbb R$, $E_1=\{\xi\in
E:\int_{-\infty}^\infty\xi(x)\,dx=1\}$, $U\in\mathcal N(E)$ and
$\gamma\in C(0)$. Let $f\in E^{(\gamma,U)}$ be arbitrary. Then every
$\zeta\in E_1$ gives $U_\zeta\in\mathcal N(E)$ and $y_\zeta\in
E^{(\gamma,U_\zeta)}$ such that $y'_\zeta=f$ and
$$y_\zeta(\xi)=f(-\int_{-\infty}^x[\xi(\tau)-(\int_{-\infty}^\infty\xi(s)\,ds)\zeta(\tau)]\,d\tau),\ \forall\,\xi\in E.$$

Proof. Only need to consider $E=\mathscr D(\mathbb R)$. Pick a
$\zeta\in E_1$ and define $A_\zeta:\mathscr D(\mathbb
R)\rightarrow\mathscr D(\mathbb R)$ by
$A_\zeta(\xi)=\xi-(\int_{-\infty}^\infty\xi(\tau)\,d\tau)\zeta,\
\forall\,\xi\in\mathscr D(\mathbb R)$. By Lemma 2.2, $A_\zeta$ is a
continuous linear operator and $\int_{-\infty}^\infty
A_\zeta(\xi)(\tau)\,d\tau=0,\ \forall\,\xi\in\mathscr D(\mathbb R)$.
For every $\xi\in\mathscr D(\mathbb R)$,
$A_\zeta(\xi)(x)=\frac{d}{dx}[\int_{-\infty}^x
A_\zeta(\xi)(\tau)\,d\tau],\ \forall\,x\in\mathbb R$. Since
$A_\zeta(\xi)\in\mathscr D(\mathbb R)$, $\forall\,\xi\in\mathscr
D(\mathbb R)$, there ia an $a_\xi>0$ such that $A_\zeta(\xi)(x)=0,\
\forall\,|x|>a_\xi$ and so
$\int_{-\infty}^xA_\zeta(\xi)(\tau)\,d\tau=0$ for $x<-a_\xi$ and
$\int_{-\infty}^xA_\zeta(\xi)(\tau)\,d\tau=\int_{-\infty}^\infty
A_\zeta(\xi)(\tau)\,d\tau=0$ for $x>a_\xi$. Thus,
$\int_{-\infty}^xA_\zeta(\xi)(\tau)\,d\tau$ gives a test function in
$\mathscr D(\mathbb R),\ \forall\,\xi\in\mathscr D(\mathbb R)$.

Let $T(\xi)(x)=\int_{-\infty}^xA_\zeta(\xi)(\tau)\,d\tau,\
\forall\,\xi\in\mathscr D(\mathbb R),\ x\in\mathbb R$. Since
$A_\zeta$ is linear, $T:\mathscr D(\mathbb R)\rightarrow\mathscr
D(\mathbb R)$ is a linear operator. Let $\xi_k\rightarrow0$ in
$\mathscr D(\mathbb R)$. By Lemma 2.2, $A_\zeta(\xi_k)\rightarrow0$
in $\mathscr D(\mathbb R)$ and so $A_\zeta(\xi_k)\rightarrow0$ in
$\mathscr D_{m_0}(\mathbb R)$ for some $m_0\in\mathbb N$ [3, p.219].
Then $\{T(\xi_k)\}\subset\mathscr D_{m_0}(\mathbb R)$. In fact,
$\{A_\zeta(\xi_k)\}\subset\mathscr D_{m_0}(\mathbb R)$, i.e.,
$A_\zeta(\xi_k)(x)=0,\ \forall\,|x|>m_0,\ k\in\mathbb N$, hence
$T(\xi_k)(x)=\int_{-\infty}^xA_\zeta(\xi_k)(\tau)\,d\tau=0$ for
$x<-m_0,\ k\in\mathbb N$ and
$T(\xi_k)(x)=\int_{-\infty}^xA_\zeta(\xi_k)(\tau)\,d\tau=\int_{-\infty}^\infty
A_\zeta(\xi_k)(\tau)\,d\tau=0$ whenever $x>m_0$ and $k\in\mathbb N$.
Since $A_\zeta(\xi_k)\rightarrow0$ in $\mathscr D_{m_0}(\mathbb R)$,
\begin{align*}
\|T(\xi_k)\|_0&=\sup_{|x|\leq m_0}|T(\xi_k)(x)| =\sup_{|x|\leq
m_0}|\int_{-m_0}^xA_\zeta(\xi_k)(\tau)\,d\tau|\\
&\leq\int_{-m_0}^{m_0}|A_\zeta(\xi_k)(\tau)|\,d\tau\leq2m_0\|A_\zeta(\xi_k)\|_0\rightarrow0,\
\ \mbox{i.e.,}\ \ \|T(\xi_k)\|_0\rightarrow0.
\end{align*}
Moreover, $\frac{dT(\xi)}{dx}(x)=A_\zeta(\xi)(x),\
\forall\,\xi\in\mathscr D(\mathbb R),\ x\in\mathbb R$, i.e.,
$\frac{dT(\xi)}{dx}=A_\zeta(\xi),\ \forall\,\xi\in\mathscr D(\mathbb
R)$. Since $\{T(\xi_k)\}\subset\mathscr D_{m_0}(\mathbb R)$ and
$A_\zeta(\xi_k)\rightarrow0$ in $\mathscr D_{m_0}(\mathbb R)$,
$\|T(\xi_k)\|_p\leq\linebreak\max\{\|T(\xi_k)\|_0,\
\|A_\zeta(\xi_k)\|_{p-1}\}\rightarrow0$ for each $p\in\mathbb N$.
Thus, $T(\xi_k)\rightarrow0$ in $\mathscr D_{m_0}(\mathbb R)$, i.e.,
$T(\xi_k)\rightarrow0$ in $\mathscr D(\mathbb R)$. Therefore,
$T:\mathscr D(\mathbb R)\rightarrow\mathscr D(\mathbb R)$ is
sequentially continuous. Since $T$ is linear and $\mathscr D(\mathbb
R)$ is bornological, i.e., $\mathscr D(\mathbb R)$ is
$C-sequential$, $T$ is continuous and so there is a balanced
$U_\zeta\in\mathcal N(\mathscr D(\mathbb R))$ such that
$T(U_\zeta)\subset U$.

Define $y_\zeta:\mathscr D(\mathbb R)\rightarrow\mathbb C$ by
$$y_\zeta(\xi)=f(-T(\xi))=
f(-\int_{-\infty}^x[\xi(\tau)-(\int_{-\infty}^\infty\xi(s)\,ds)\zeta(\tau)]d\tau),
\ \forall\,\xi\in\mathscr D(\mathbb R).$$ Since both $f$ and $T$ are
continuous, $y_\zeta$ is continuous.

Let $\xi\in\mathscr D(\mathbb R)$, $\eta\in U_\zeta$ and $|t|\leq1$.
Since $U_\zeta$ is balanced and $f\in\mathscr D(\mathbb
R)^{(\gamma,U)}$, $T(-\eta)\in T(U_\zeta)\subset U$ and
\begin{align*}
y_\zeta(\xi+t\eta)&=f(-T(\xi+t\eta))=f(-T(\xi)+tT(-\eta))=rf(-T(\xi))+sf(T(-\eta))\\
&=rf(-T(\xi))+sf(-T(\eta))=ry_\zeta(\xi)+sy_\zeta(\eta),\
|r-1|\leq|\gamma(t)|,\,|s|\leq|\gamma(t)|.
\end{align*}
Thus, $y_\zeta\in\mathscr D(\mathbb R)^{(\gamma,U_\zeta)}$.

For every $\xi\in\mathscr D(\mathbb R)$,
$T(\xi')(x)=\int_{-\infty}^xA_\zeta(\xi')(\tau)\,d\tau
=\int_{-\infty}^x[\xi'(\tau)-(\int_{-\infty}^\infty\xi'(s)\,ds)\zeta(\tau)]\linebreak
d\tau =\int_{-\infty}^x\xi'(\tau)\,d\tau=\xi(x),\
\forall\,x\in\mathbb R$, i.e., $T(\xi')=\xi,\
\forall\,\xi\in\mathscr D(\mathbb R)$. Then
$$y'_\zeta(\xi)=y_\zeta(-\xi')=f(-T(-\xi'))=f(T(\xi'))=f(\xi),\ \forall\,\xi\in\mathscr D(\mathbb R),$$
i.e., $y'_\zeta=f$. $\square$
\end{Theorem}

\begin{Corollary}
Let $E\in\{\mathscr D_a(\mathbb R),\mathscr D(\mathbb R)\}$ and
$f\in E^{(\gamma,U)}$. Then every $\zeta\in E$ which $\zeta\neq\xi'$
for all $\xi\in E$ (i.e., $\int_{-\infty}^\infty\zeta(x)\,dx\neq0$)
gives a balanced $U_\zeta\in\mathcal N(E)$ and a solution $y_\zeta$
of the equation $y'=f$ such that $y_\zeta\in E^{(\gamma,U_\zeta)}$
and
$$y_\zeta(\xi)=f(-\int_{-\infty}^x[\xi(\tau)
-\frac{\int_{-\infty}^\infty\xi(s)\,ds}{\int_{-\infty}^\infty\zeta(s)\,ds}\zeta(\tau)]\,d\tau),\
\forall\,\xi\in E.$$
\end{Corollary}

\begin{Theorem}
Let $E\in\{\mathscr D_a(\mathbb R),\mathscr D(\mathbb R)\}$ and
$f\in E^{(\gamma,U)}$. For every $\xi_0\in E$, $f_0\in
E^{(\gamma,U)}$ and $\zeta\in E$ with
$\int_{-\infty}^\infty\zeta(x)\,dx\neq0$ let
$$g(\xi)=f_0[(\int_{-\infty}^\infty\xi(x)\,dx)\xi_0]+f(-\int_{-\infty}^x[\xi(\tau)
-\frac{\int_{-\infty}^\infty\xi(s)\,ds}{\int_{-\infty}^\infty\zeta(s)\,ds}\zeta(\tau)]\,d\tau),\
\forall\,\xi\in E,$$ then $g\in[E^{(\gamma,W)}]$ for some
$W\in\mathcal N(E)$ and $g'=f$.

Proof. Let $\xi_0\in E$, $f_0\in E^{(\gamma,U)}$ and $\zeta\in E$
with $\int_{-\infty}^\infty\zeta(x)\,dx\neq0$. Let
$$g_0(\xi)=f_0[(\int_{-\infty}^\infty\xi(x)\,dx)\xi_0],\ \xi\in E,$$
$$y_\zeta(\xi)=f(-\int_{-\infty}^x[\xi(\tau)
-\frac{\int_{-\infty}^\infty\xi(s)\,ds}{\int_{-\infty}^\infty\zeta(s)\,ds}\zeta(\tau)]\,d\tau),\
\xi\in E.$$ Then $g_0\in E^{(\gamma,V)}$ for some $V\in\mathcal
N(E)$ and $g'_0=0$ by Th. 2.3, and $y_\zeta\in E^{(\gamma,U_\zeta)}$
for some balanced $U_\zeta\in\mathcal N(E)$ and $y'_\zeta=f$ by Cor.
2.3.

Let $W=V\bigcap U_\zeta$. Then $W\in\mathcal N(E)$ and
$E^{(\gamma,V)}\bigcup E^{(\gamma,U_\zeta)}\subset E^{(\gamma,W)}$.
Thus, $g_0\in E^{(\gamma,W)}$, $y_\zeta\in E^{(\gamma,W)}$,
$g=g_0+y_\zeta\in[E^{(\gamma,W)}]$ and
$g'=(g_0+y_\zeta)'=g'_0+y'_\zeta=f$. $\square$
\end{Theorem}

Further discussions of ordinary and partial differential equations
of demi-distributions will be interesting but we reserve these
discussions for another paper.

\section{Fourier Transform}
Let $x+iy=(x_1+iy_1,\cdots,x_n+iy_n)\in\mathbb C^n$,
$|y|=|y_1|+\cdots+|y_n|$. For $a>0$ and
$\alpha=(\alpha_1,\cdots,\alpha_n)$, a multi-index, let
$(x+iy)^\alpha=\prod_{k=1}^n(x_k+iy_k)^{\alpha_k}$ and
$$Z(a)=\big\{\zeta\in\mathbb C^{\mathbb C^n}:\zeta\mbox{ is entire; for every multi-index }\alpha,\
|(x+iy)^\alpha\zeta(x+iy)|\leq C_\alpha(\zeta)e^{a|y|}\big\},$$
$$\|\zeta\|_p=\sup_{x+iy\in\mathbb C^n,|\alpha|\leq p}|(x+iy)^\alpha\zeta(x+iy)|e^{-a|y|},\ p=0,1,2,3,\cdots,$$
$$Z=\big\{\zeta\in\mathbb C^{\mathbb C^n}:\zeta\mbox{ is entire, }\exists\,a(\zeta)>0
\mbox{ such that }|(x+iy)^\alpha\zeta(x+iy)|\leq
C_\alpha(\zeta)e^{a(\zeta)|y|}\big\}.$$

The Fourier transform $F(\xi)$ of $\xi\in\mathscr D$ is given by
$$F(\xi)(\sigma+i\tau)=\zeta(\sigma+i\tau)=\int\xi(x)e^{i(x,\sigma)-(x,\tau )}\,dx,\ (x,\sigma)
=\sum_{j=1}^nx_j\sigma_j,\ (x,\tau)=\sum_{j=1}^nx_j\tau_j.$$ Then
$F[\mathscr D_a]=Z_a$, $F[\mathscr D]=Z$, and operators
\begin{align*}
F:\mathscr D_a\rightarrow Z_a,\ \ F:\mathscr D\rightarrow Z,\ \
F^{-1}:Z_a\rightarrow\mathscr D_a,\ \ F^{-1}:Z\rightarrow\mathscr D
\end{align*}
are both continuous and linear [4, 3.1.1---3.1.2].

For $\xi\in\mathscr S$ let
$$F(\xi)(\sigma)=\zeta(\sigma)=\int\xi(x)e^{i(x,\sigma)}\,dx,\ \forall\,\sigma\in\mathbb R^n.$$
Then $F(\mathscr S)=\mathscr S$ and both $F$ and $F^{-1}$ are
continuous and linear.

\begin{Definition}
Let $E\in\{\mathscr D_a,\mathscr D,\mathscr S\}$, $U\in\mathcal
N(E)$ and $\gamma\in C(0)$. For $f\in E^{(\gamma,U)}$ define
$\hat{f}:F(E)\rightarrow\mathbb C$ by
$\hat{f}(\zeta)=(2\pi)^nf(F^{-1}(\zeta)),\ \forall\,\zeta\in F(E)$.
We write $\hat{f}=F(f)$ and so
$$F(f)(F(\xi))=(2\pi)^nf(\xi),\ \forall\,f\in E^{(\gamma,U)},\ \xi\in E.$$
\end{Definition}

Henceforth, $E\in\{\mathscr D_a,\mathscr D,\mathscr S\}$,
$U\in\mathcal N(E)$ and $\gamma\in C(0)$.

\begin{Theorem}
$F(E^{(\gamma,U)})=(F(E))^{(\gamma,F(U))}.$

Proof. Since both $F$ and $F^{-1}$ are continuous linear operators,
$F(U)\in\mathcal N(F(E))$. Let $f\in E^{(\gamma,U)}$, $\zeta\in
F(E)$, $\eta\in F(U)$ and $|t|\leq1$. Then
\begin{align*}
F(f)(\zeta+t\eta)&=(2\pi)^nf(F^{-1}(\zeta+t\eta))=(2\pi)^nf(F^{-1}(\zeta)+tF^{-1}(\eta))\\
&=(2\pi)^n[rf(F^{-1}(\zeta))+sf(F^{-1}(\eta))]\\
&=r(2\pi)^nf(F^{-1}(\zeta))+s(2\pi)^nf(F^{-1}(\eta))\\
&=rF(f)(\zeta)+sF(f)(\eta),\qquad|r-1|\leq|\gamma(t)|,\
|s|\leq|\gamma(t)|.
\end{align*}
Thus, $F(f)\in\mathscr L_{\gamma,F(U)}(F(E),\mathbb C)$.

Let $\zeta_\alpha\rightarrow\zeta$ in $F(E)$. Then
$F^{-1}(\zeta_\alpha)\rightarrow F^{-1}(\zeta)$ in $E$ and so
$F(f)(\zeta_\alpha)=(2\pi)^nf(F^{-1}(\zeta_\alpha))\rightarrow(2\pi)^nf(F^{-1}(\zeta))=F(f)(\zeta)$.
This shows that $F(f)$ is continuous and $F(f)\in
(F(E))^{(\gamma,F(U))}$.

Conversely, for $g\in(F(E))^{(\gamma,F(U))}$ define
$$f(\xi)=(2\pi)^{-n}g(F(\xi)),\
\forall\,\xi\in E.$$ If $\xi\in E,\ \eta\in U$ and $|t|\leq1$, then
\begin{align*}
f(\xi+t\eta)&=(2\pi)^{-n}g(F(\xi+t\eta))=(2\pi)^{-n}g(F(\xi)+tF(\eta))\\
&=(2\pi)^{-n}[rg(F(\xi))+sg(F(\eta))]\\
&=rf(\xi)+sf(\eta),\qquad|r-1|\leq|\gamma(t)|,\ |s|\leq|\gamma(t)|,
\end{align*}
i.e., $f\in\mathscr L_{\gamma,U}(E,U)$.

Since both $g$ and $F$ are continuous, $f$ is continuous so $f\in
E^{(\gamma,U)}$ and $F(f)(\zeta)=(2\pi)^nf(F^{-1}(\zeta))=g(\zeta),\
\forall\,\zeta\in F(E)$, i.e., $g=F(f)$. $\square$
\end{Theorem}

\begin{Definition}
Let $[(F(E))^{(\gamma,F(U))}]=span\,(F(E))^{(\gamma,F(U))}$ in
$\mathbb C^{F(E)}$. For $f\in[E^{(\gamma,U)}]$ define
$F(f):F(E)\rightarrow\mathbb C$ by
$$F(f)(F(\xi))=(2\pi)^nf(\xi),\ \forall\,\xi\in E.$$
\end{Definition}

\begin{Theorem}
If $f=\sum_{k=1}^m\alpha_kf_k$ where $\alpha_k\in\mathbb C$ and
$f_k\in E^{(\gamma,E)}$, then
$F(f)=\sum_{k=1}^m\alpha_kF(f_k)\in[(F(E))^{(\gamma,F(U))}]$, and
$$F([E^{(\gamma,U)}])=[(F(E))^{(\gamma,F(U))}].$$
Moreover, $F:[E^{(\gamma,U)}]\rightarrow[(F(E))^{(\gamma,F(U))}]$ is
$w\ast-w\ast$ continuous and linear.

Proof. For $\xi\in E$,
$F(f)(F(\xi))=(2\pi)^nf(\xi)=(2\pi)^n\sum_{k=1}^m\alpha_kf_k(\xi)
=\linebreak\sum_{k=1}^m\alpha_kF(f_k)(F(\xi))=(\sum_{k=1}^m\alpha_kF(f_k))(F(\xi))$.
Thus, $F(f)=\sum_{k=1}^m\alpha_kF(f_k)\in [(F(E))^{(\gamma,F(U))}]$.

Let $g=\sum_{k=1}^m\alpha_kg_k\in[(F(E))^{(\gamma,F(U))}]$ where
$\alpha_k\in\mathbb C$, $g_k\in(F(E))^{(\gamma,F(U))}$. By Th. 3.1,
each $g_k=F(f_k)$ for some $f_k\in E^{(\gamma,U)}$ and so
$g=\sum_{k=1}^m\alpha_kg_k=\sum_{k=1}^m\alpha_kF(f_k)=F(\sum_{k=1}^m\alpha_kf_k)\in
F([E^{(\gamma,U)}])$ because
$(\sum_{k=1}^m\alpha_kF(f_k))(F(\xi))=\linebreak\sum_{k=1}^m\alpha_kF(f_k)(F(\xi))
=(2\pi)^n\sum_{k=1}^m\alpha_kf_k(\xi)=(2\pi)^n(\sum_{k=1}^m\alpha_kf_k)(\xi)=\linebreak
F(\sum_{k=1}^m\alpha_kf_k)(F(\xi)),\ \forall\,\xi\in E$.

If $f_\alpha\stackrel{w\ast}{\longrightarrow}f$ in
$[E^{(\gamma,U)}]$, i.e., $f_\alpha(\xi)\rightarrow f(\xi),\
\forall\,\xi\in E$, then
$F(f_\alpha)(F(\xi))=(2\pi)^nf_\alpha(\xi)\rightarrow(2\pi)^nf(\xi)=F(f)(F(\xi)),\
\forall\,\xi\in E$, i.e.,
$F(f_\alpha)\stackrel{w\ast}{\longrightarrow}F(f)$.

For $f,h\in[E^{(\gamma,U)}]$ and $\alpha,\beta\in\mathbb C$,
\begin{align*}
F(\alpha f+\beta h)(F(\xi))&=(2\pi)^n(\alpha f+\beta h)(\xi)\\
&=(2\pi)^n(\alpha f(\xi)+\beta h(\xi))\\
&=\alpha F(f)(F(\xi))+\beta F(h)(F(\xi))\\
&=(\alpha F(f)+\beta F(h))(F(\xi)),\ \forall\,\xi\in E,
\end{align*}
i.e., $F(\alpha f+\beta h)=\alpha F(f)+\beta F(h)$. $\square$
\end{Theorem}

Now we consider the case of $n=1$. Let $\mathscr S(\mathbb R)$ be
the space of infinitely differentiable but rapidly decreasing
functions defined on $\mathbb R$. Then $F(\mathscr S(\mathbb
R))=\mathscr S(\mathbb R)$, $F((\mathscr S(\mathbb R))')=(\mathscr
S(\mathbb R))'$.

A constant $C\in(\mathscr S(\mathbb R))'$ means that
$C(\zeta)=\int_{-\infty}^\infty C\zeta(\sigma)\,d\sigma$ for all
$\zeta\in\mathscr S(\mathbb R)$ [4, 3.2.1], and
$C=F(C\delta)=CF(\delta)$. In fact, $CF(\delta)(F(\xi))=2\pi
C\delta(\xi)=2\pi C\xi(0)=2\pi C\frac{1}{2\pi}\int_{-\infty}^\infty
e^{-i0\sigma}F(\xi)(\sigma)\,d\sigma=\int_{-\infty}^\infty
CF(\xi)(\sigma)\,d\sigma=C(F(\xi))$ for all $\xi\in\mathscr
S(\mathbb R)$.

\begin{Lemma}
Let $y\in(\mathscr S(\mathbb R))'$, a usual distribution. Then
$$y(i\sigma\zeta(\sigma))=0\mbox{ for all }\zeta\in\mathscr S(\mathbb R)$$
if and only if $y=C\delta$, where $C$ is a constant.

Proof. Suppose that $y\in(\mathscr S(\mathbb R))'$ and
$y(i\sigma\zeta(\sigma))=0,\ \forall\,\zeta\in\mathscr S(\mathbb
R)$. Since $(\mathscr S(\mathbb R))'=F((\mathscr S(\mathbb R))')$,
there is a usual distribution $f\in(\mathscr S(\mathbb R))'$ such
that $y=F(f)$ and
\begin{align*}
f'(\xi)=f(-\xi')&=\frac{1}{2\pi}F(f)(F((-\xi)'))=\frac{1}{2\pi}F(f)(-i\sigma
F(-\xi)(\sigma))\\
&=\frac{1}{2\pi}y(i\sigma F(\xi)(\sigma))=0,\
\forall\,\xi\in\mathscr S(\mathbb R),
\end{align*}
i.e., $f'=0$. But $f$ is a usual distribution so $f=C$, a constant.
Thus, $y=F(f)=F(C)=CF(1)=C\delta$.

Conversely, if $y=C\delta$ where $C$ is a constant, then
$$y(i\sigma\zeta(\sigma))=C\delta(i\sigma\zeta(\sigma))=0,\ \forall\,\zeta\in\mathscr S(\mathbb R).\ \square$$
\end{Lemma}

However, there exists a lot of various demi-distributions on
$\mathscr S(\mathbb R)$ which satisfy the condition
$y(i\sigma\zeta(\sigma))=0,\ \forall\,\zeta\in\mathscr S(\mathbb
R)$.

\begin{Theorem}
Let $U\in\mathcal N(\mathscr S(\mathbb R))$ and $\gamma\in C(0)$.
Pick an arbitrary $f_0\in(\mathscr S(\mathbb R))^{(\gamma,U)}$ and
$\xi_0\in\mathscr S(\mathbb R)$ and let
$f(\xi)=f_0[(\int_{-\infty}^\infty\xi(t)\,dt)\xi_0],\
\forall\,\xi\in\mathscr S(\mathbb R)$. Then $f\in(\mathscr S(\mathbb
R))^{(\gamma,V)}$ for some $V\in\mathcal N(\mathscr S(\mathbb R))$
and $F(f)\in(\mathscr S(\mathbb R))^{(\gamma,F(V))}$ such that
$$F(f)(i\sigma\zeta(\sigma))=0,\ \forall\,\zeta\in\mathscr S(\mathbb R).$$

Proof. By Th. 2.3, there is a $V\in\mathcal N(\mathscr S(\mathbb
R))$ such that $f\in(\mathscr S(\mathbb R))^{(\gamma,V)}$ and
$f'=0$. If $\zeta\in\mathscr S(\mathbb R)$, then $\zeta=F(\xi)$ for
some $\xi\in\mathscr S(\mathbb R)$ and
\begin{align*}
F(f)(i\sigma\zeta(\sigma))&=F(f)(i\sigma
F(\xi)(\sigma))=F(f)(-i\sigma F(-\xi)(\sigma))\\
&=F(f)(F((-\xi)'))=2\pi f(-\xi')=2\pi f'(\xi)=0.\ \square
\end{align*}
\end{Theorem}

\begin{Example}
Let $\gamma(t)=\frac{\pi}{2}\,t$ for $t\in\mathbb R$,
$U=\{\eta\in\mathscr S(\mathbb R):|\eta(0)|<1\}$. Then
$f_0=\sin\,\circ\,\delta\in(\mathscr S(\mathbb R))^{(\gamma,U)}$,
where $f_0(\xi)=\sin\,[\delta(\xi)]=\sin\,\xi(0)$ for all
$\xi\in\mathscr S(\mathbb R)$. Let
$\xi_0=\begin{cases}e^{\frac{1}{x^2-1}},&|x|<1,\\0,&|x|\geq1\end{cases}$
and $f(\xi)=f_0[(\int_{-\infty}^\infty\xi(t)\,dt)\xi_0]$ for
$\xi\in\mathscr S(\mathbb R)$, i.e.,
$f(\xi)=\sin\,[\delta((\int_{-\infty}^\infty\xi(t)\,dt)\xi_0)]=\sin\,(e^{-1}\int_{-\infty}^\infty\xi(t)\,dt)$,
$\forall\,\xi\in\mathscr S(\mathbb R)$. Pick a $V\in\mathcal
N(\mathscr S(\mathbb R))$ for which
$(\int_{-\infty}^\infty\eta(t)\,dt)\xi_0\in U,\ \forall\,\eta\in V$.
Then $f\in(\mathscr S(\mathbb R))^{(\gamma,V)}$ and
$f'(\xi)=f(-\xi')=\sin\,(-e^{-1}\int_{-\infty}^\infty\xi'(t)\,dt)=0$
for all $\xi\in\mathscr S(\mathbb R)$, i.e., $f'=0$. Therefore,
$F(f)(i\sigma\zeta(\sigma))=0,\ \forall\,\zeta\in\mathscr S(\mathbb
R)$.

By Th. 3.1, $F(f)\in(\mathscr S(\mathbb R))^{(\gamma,F(V))}$ and
$$F(f)(\zeta)=2\pi f(F^{-1}(\zeta))
=2\pi \sin\,[e^{-1}\int_{-\infty}^\infty\int_{-\infty}^\infty
e^{-ix\sigma}\zeta(\sigma)\,d\sigma\,dx], \
\forall\,\zeta\in\mathscr S(\mathbb R),$$
$$F(f)(F(\xi))=2\pi f(\xi)
=2\pi \sin\,[e^{-1}\int_{-\infty}^\infty\xi(x)\,dx], \
\forall\,\xi\in\mathscr S(\mathbb R).$$ Since
$C\delta(F(\xi))=C\delta(\int_{-\infty}^\infty
e^{ix\sigma}\xi(x)\,dx)=C\int_{-\infty}^\infty\xi(x)\,dx$ for all
$\xi\in\mathscr S(\mathbb R)$, $F(f)\neq C\delta$ for every constant
$C$.
\end{Example}

\section{Convolutions}
In this section, $E\in\{\mathscr D,\mathscr S\}$, $U\in\mathcal
N(E)$ and $\gamma\in C(0)$.

\begin{Definition}
([4, 3.3.2]) A distribution $f_0\in E'$ is called a convolution
multiplier on $E$ if the following (i) and (ii) hold for $f_0$:

(i) if for each $\xi\in E$ define $f_0\ast\xi:\mathbb
R^n\rightarrow\mathbb C$ by
$$(f_0\ast\xi)(x)=f_0(\xi(x+\cdot)),\
\forall\,x\in\mathbb R^n,$$ then $f_0\ast\xi\in E$;

(ii) if $\xi_k\rightarrow0$ in $E$, then $f_0\ast\xi_k\rightarrow0$
in $E$.
\end{Definition}

\begin{Lemma}
If $f_0$ is a convolution multiplier on $E$, then
$f_0\ast\cdot:E\rightarrow E$ is a continuous linear operator.

Proof. For $\xi,\eta\in E$ and $t\in\mathbb C$,
$(f_0\ast(\xi+t\eta))(x)=f_0(\xi(x+\cdot)+t\eta(x+\cdot))=f_0(\xi(x+\cdot))+tf_0(\eta(x+\cdot))
=(f_0\ast\xi)(x)+t(f_0\ast\eta)(x),\ \forall\,x\in\mathbb R^n$,
i.e., $f_0\ast(\xi+t\eta)=f_0\ast\xi+t(f_0\ast\eta)$.

By Def. 4.1, $f_0\ast\cdot$ is sequentially continuous. Since
$\mathscr S$ is metric and $\mathscr D$ is bornological,
$f_0\ast\cdot$ is continuous. $\square$
\end{Lemma}

Following [4], $P(D)=\sum a_\alpha D^\alpha=\sum
a_{\alpha_1,\cdots,\alpha_n}\frac{\partial^{\alpha_1+\cdots+\alpha_n}}{\partial
x_1^{\alpha_1}\cdots\partial x_n^{\alpha_n}}$ is a finite sum, where
$a_\alpha\in\mathbb C$ and $\alpha=(\alpha_1,\cdots,\alpha_n)$ is a
multi-index. For $\xi\in E$ and $x\in\mathbb R^n$,
$\frac{\partial\xi(x+\tau)}{\partial\tau_j}
=\frac{\partial\xi(x+\tau)}{\partial(x_j+\tau_j)}\frac{\partial(x_j+\tau_j)}{\partial\tau_j}
=\frac{\partial\xi(x+\tau)}{\partial(x_j+\tau_j)}(\frac{\partial
x_j}{\partial\tau_j}+\frac{\partial\tau_j}{\partial\tau_j})=\frac{\partial\xi(x+\tau)}{\partial(x_j+\tau_j)}$
and, by induction, it is easy to see that
$\frac{\partial^{|\alpha|}\xi(x+\tau)}{\partial\tau_1^{\alpha_1}\cdots\partial\tau_n^{\alpha_n}}
=\frac{\partial^{|\alpha|}\xi(x+\tau)}{\partial(x_1+\tau_1)^{\alpha_1}\cdots\partial(x_n+\tau_n)^{\alpha_n}}$
for every multi-index $\alpha$ and so there is no any ambiguity for
the notation $D^\alpha\xi(x+\cdot)$, i.e.,
$$D^\alpha\xi(x+\cdot)=\frac{\partial^{|\alpha|}\xi(x+\tau)}{\partial\tau_1^{\alpha_1}\cdots\partial\tau_n^{\alpha_n}}
=\frac{\partial^{|\alpha|}\xi(x+\tau)}{\partial(x_1+\tau_1)^{\alpha_1}\cdots\partial(x_n+\tau_n)^{\alpha_n}}
=(D^\alpha\xi)(x+\cdot).$$

\begin{Theorem}
If $f_0\in E'$ is a convolution multiplier on $E$, then $P(D)f_0$ is
also a convolution multiplier on $E$, and
$$(\sum a_\alpha D^\alpha f_0)\ast\xi=\sum a_\alpha(-1)^{|\alpha|}f_0\ast D^\alpha\xi,\ \forall\,\xi\in E.$$

Proof. For $\xi\in E$, define $(\sum a_\alpha D^\alpha
f_0)\ast\xi:\mathbb R^n\rightarrow\mathbb C$ by
$$((\sum a_\alpha D^\alpha f_0)\ast\xi)(x)=(\sum a_\alpha D^\alpha f_0)(\xi(x+\cdot))
=\sum a_\alpha(D^\alpha f_0)(\xi(x+\cdot)),\ \forall\,x\in\mathbb
R^n.$$ For $x\in\mathbb R^n$ we write $\xi(x+\tau)=\zeta_x(\tau),\
\forall\,\tau\in\mathbb R^n$. Then $\zeta_x\in E$ and $(D^\alpha
f_0)(\xi(x+\cdot))=(D^\alpha
f_0)(\zeta_x)=f_0((-1)^{|\alpha|}D^\alpha\zeta_x)
=f_0((-1)^{|\alpha|}\frac{\partial^{|\alpha|}\zeta_x}{\partial\tau_1^{\alpha_1}\cdots\partial\tau_n^{\alpha_n}})
=f_0((-1)^{|\alpha|}\frac{\partial^{|\alpha|}\xi(x+\tau)}{\partial\tau_1^{\alpha_1}\cdots\partial\tau_n^{\alpha_n}})
\linebreak=f_0[(-1)^{|\alpha|}(D^\alpha\xi)(x+\cdot)]$,
\begin{align*}
((\sum a_\alpha D^\alpha f_0)\ast\xi)(x)&=\sum a_\alpha(D^\alpha
f_0)(\xi(x+\cdot))\\
&=\sum a_\alpha f_0[(-1)^{|\alpha|}(D^\alpha\xi)(x+\cdot)]\\
&=\sum a_\alpha(-1)^{|\alpha|}f_0[(D^\alpha\xi)(x+\cdot)]\\
&=\sum a_\alpha(-1)^{|\alpha|}(f_0\ast D^{\alpha}\xi)(x)\\
&=(\sum a_\alpha(-1)^{|\alpha|}f_0\ast D^\alpha\xi)(x),\
\forall\,x\in\mathbb R^n.
\end{align*}
Thus, $(\sum a_\alpha D^\alpha f_0)\ast\xi=\sum
a_\alpha(-1)^{|\alpha|}f_0\ast D^\alpha\xi\in E,\ \forall\,\xi\in
E$.

Let $\xi_k\rightarrow0$ in $E$. By Lemma 2.1,
$D^\alpha\xi_k\rightarrow0$ in $E$ and
$$\lim_k(\sum a_\alpha D^\alpha f_0)\ast\xi_k=\sum a_\alpha(-1)^{|\alpha|}\lim_k(f_0\ast D^\alpha\xi_k)=0.\ \square$$
\end{Theorem}

\begin{Definition}
For every convolution multiplier $f_0\in E'$ and
$f\in[E^{(\gamma,U)}]$ define the convolution $f_0\ast
f:E\rightarrow\mathbb C$ by
$$(f_0\ast f)(\xi)=f(f_0\ast\xi),\ \forall\,\xi\in E.$$
\end{Definition}

\begin{Example}
(1) For every $\xi\in E$, $\delta\ast\xi=\xi$,
$(D^\alpha\delta)\ast\xi=(-1)^{|\alpha|}D^\alpha\xi$, and for every
$f\in[E^{(\gamma,U)}]$, $\delta\ast f=f$, $(D^\alpha\delta)\ast
f=D^\alpha f$.

In fact, for $\xi\in E$ and $x\in\mathbb R^n$,
$(\delta\ast\xi)(x)=\delta(\xi(x+\cdot))=\xi(x+0)=\xi(x)$,
$((D^\alpha\delta)\ast\xi)(x)=(D^\alpha\delta)(\xi(x+\cdot))
=\delta((-1)^{|\alpha|}\frac{\partial^{|\alpha|}\xi(x+\tau)}{\partial\tau_1^{\alpha_1}\cdots\partial\tau_n^{\alpha_n}})
=(-1)^{|\alpha|}(D^\alpha\xi)(x)$, i.e., $\delta\ast\xi=\xi$,
$(D^\alpha\delta)\ast\xi=(-1)^{|\alpha|}D^\alpha\xi$. Then for every
$\xi\in E$,
$$(\delta\ast f)(\xi)=f(\delta\ast\xi)=f(\xi),\ ((D^\alpha\delta)\ast f)(\xi)=f((D^\alpha\delta)\ast\xi)
=f((-1)^{|\alpha|}D^\alpha\xi)=(D^\alpha f)(\xi).$$

(2) Let $f\in E'$ and $U=\{\eta\in E:|f(\eta)|<1\}$. By Cor. 1.1,
$h\circ f\in E^{(\gamma,U)}$ for each $h\in\mathscr
L_{\gamma,1}(\mathbb C,\mathbb C)$ and $f_0\ast(h\circ
f)=h\circ(f_0\ast f)$ for every convolution multiplier $f_0\in E'$.
In fact,
$$(f_0\ast(h\circ f))(\xi)=(h\circ f)(f_0\ast\xi)=h[f(f_0\ast\xi)]
=h[(f_0\ast f)(\xi)]=[h\circ(f_0\ast f)](\xi),\ \forall\,\xi\in E.$$
\end{Example}

\begin{Theorem}
Let $f_0\in E'$ be a convolution multiplier. For every $U\in\mathcal
N(E)$ there is a $V\in\mathcal N(E)$ such that
$$f_0\ast f\in E^{(\gamma,V)},\ \forall\,f\in E^{(\gamma,U)},$$
$$f_0\ast f\in[E^{(\gamma,V)}],\ \forall\,f\in[E^{(\gamma,U)}].$$

Proof. Let $T(\xi)=f_0\ast\xi,\ \forall\,\xi\in E$. By Lemma 4.1,
$T$ is a continuous linear operator. Let $U\in\mathcal N(E)$ and
$f=\sum_{k=1}^m\alpha_k f_k\in[E^{(\gamma,U)}]$ where
$\alpha_k\in\mathbb C$ and $f_k\in E^{(\gamma,U)}$. There is a
$V\in\mathcal N(E)$ such that $T(V)\subset U$, i.e., $f_0\ast\eta\in
U,\ \forall\,\eta\in V$. Since $(f_0\ast
f)(\xi)=f(f_0\ast\xi)=f(T(\xi))$ for all $\xi\in E$ and both $f$ and
$T$ are continuous, $f_0\ast f:E\rightarrow\mathbb C$ is continuous.

Let $\xi\in E$, $\eta\in V$ and $|t|\leq1$. For $1\leq k\leq m$,
\begin{align*}
(f_0\ast f_k)(\xi+t\eta)&=f_k[f_0\ast(\xi+t\eta)]=f_k[T(\xi+t\eta)]\\
&=f_k(T(\xi)+tT(\eta))=r_kf_k(T(\xi))+s_kf_k(T(\eta))\\
&=r_k(f_0\ast f_k)(\xi)+s_k(f_0\ast f_k)(\eta),
\end{align*}
where $|r_k-1|\leq|\gamma(t)|$, $|s_k|\leq|\gamma(t)|$. Thus,
$f_0\ast f_k\in E^{(\gamma, V)},\ k=1,2,\cdots,m$.

Now $(f_0\ast f)(\xi)=f(f_0\ast\xi)=(\sum_{k=1}^m\alpha_k
f_k)(f_0\ast\xi)=\sum_{k=1}^m\alpha_kf_k(f_0\ast\xi)=\sum_{k=1}^m\alpha_k(f_0\ast
f_k)(\xi)=(\sum_{k=1}^m\alpha_k(f_0\ast f_k))(\xi),\ \forall\,\xi\in
E$. This shows that $f_0\ast f=\sum_{k=1}^m\alpha_k(f_0\ast
f_k)\in[E^{(\gamma,V)}]$. $\square$
\end{Theorem}

Henceforth, we write $f_0\ast\xi=f_0(\xi(x+\cdot))=(f_0\ast\xi)(x)$,
see [4, 3.3.2].

Observe that $tf\in E^{(\gamma,U)}$ whenever $t\in\mathbb C$ and
$f\in E^{(\gamma,U)}$.

\begin{Lemma}
If $f_0\in E'$ is a convolution multiplier and $t\in\mathbb C$, then
$$t(f_0\ast\xi)=(tf_0)\ast\xi,\ \forall\,\xi\in E;$$
$$t(f_0\ast f)=f_0\ast(tf),\ \forall\,f\in[E^{(\gamma,U)}];$$
$$t(f_0\ast f)=(tf_0)\ast f,\ \forall\,f\in E'.$$

Proof.
$t(f_0\ast\xi)=tf_0(\xi(x+\cdot))=(tf_0)(\xi(x+\cdot))=(tf_0)\ast\xi,\
\forall\,\xi\in E$. For $f\in[E^{(\gamma,U)}]$ and $\xi\in E$,
$t(f_0\ast
f)(\xi)=tf(f_0\ast\xi)=(tf)(f_0\ast\xi)=(f_0\ast(tf))(\xi)$. If
$f\in E'$, then $t(f_0\ast
f)(\xi)=tf(f_0\ast\xi)=f(t(f_0\ast\xi))=f((tf_0)\ast\xi)=((tf_0)\ast
f)(\xi)$ for all $\xi\in E$. $\square$
\end{Lemma}

As usual, $e_1=(1,0,\cdots,0)$, $e_2=(0,1,0,\cdots,0)$, $\cdots$,
$e_n=(0,\cdots,0,1)$.

\begin{Theorem}
If $f_0\in E'$ is a convolution multiplier and $\alpha$ is a
multi-index, then $D^\alpha(f_0\ast\xi)=f_0\ast D^\alpha\xi$ for
$\xi\in E$, and
$$D^\alpha(f_0\ast f)=(D^\alpha f_0)\ast f=f_0\ast D^\alpha f,\ \forall\,f\in[E^{(\gamma,U)}].$$

Proof. As in [4], for $\xi\in E$ and $1\leq j\leq n$ we have
$x+he_j=(x_1,\cdots,x_n)+(0,\cdots,0,\stackrel{(j)}{h},0,\cdots,0)$
and $\frac{\partial(f_0\ast\xi)}{\partial
x_j}=\lim_{h\rightarrow0}\frac{1}{h}[(f_0\ast\xi)(x+he_j)-(f_0\ast
f)(x)]=\lim_{h\rightarrow0}\frac{1}{h}[f_0(\xi(x+he_j+\cdot))-f_0(\xi(x+\cdot))]
=\lim_{h\rightarrow0}f_0(\frac{\xi(x+he_j+\cdot)-\xi(x+\cdot)}{h})
=\linebreak
f_0(\lim_{h\rightarrow0}\frac{\xi(x+(\cdot+he_j))-\xi(x+\cdot)}{h})
=f_0(\frac{\partial\xi}{\partial\tau_j}(x+\tau))=f_0((D^{e_j}\xi)(x+\cdot))=f_0\ast
D^{e_j}\xi=f_0\ast\frac{\partial\xi}{\partial x_j}$, i.e.,
$\frac{\partial(f_0\ast\xi)}{\partial
x_j}=f_0\ast\frac{\partial\xi}{\partial x_j}$ [4, Th. 3.3.3]. If
$D^\alpha(f_0\ast\xi)=f_0\ast D^\alpha\xi$, then
$D^{\alpha+e_j}(f_0\ast\xi)=\frac{\partial}{\partial
x_j}(D^\alpha(f_0\ast\xi))=\frac{\partial}{\partial x_j}(f_0\ast
D^\alpha\xi)=f_0\ast\frac{\partial D^\alpha\xi}{\partial
x_j}=f_0\ast D^{\alpha+e_j}\xi$. Thus, $D^\alpha(f_0\ast\xi)=f_0\ast
D^\alpha\xi$ for every multi-index $\alpha$.

Let $f\in[E^{(\gamma,U)}]$, $\xi\in E$ and $1\leq j\leq n$. By Th.
4.2, $f_0\ast f\in[E^{(\gamma,V)}]$ for some $V\in\mathcal N(E)$ and
$\frac{\partial(f_0\ast f)}{\partial x_j}(\xi)=(f_0\ast
f)(-\frac{\partial\xi}{\partial
x_j})=f(f_0\ast(-1)\frac{\partial\xi}{\partial
x_j})=f[f_0(-\frac{\partial\xi}{\partial
\tau_j}(x+\tau))]=f[\frac{\partial
f_0}{\partial\tau_j}(\xi(x+\tau))]=f((D^{e_j}f_0)\ast\xi)=((D^{e_j}f_0)\ast
f)(\xi)=(\frac{\partial f_0}{\partial x_j}\ast f)(\xi)$. Thus,
$\frac{\partial(f_0\ast f)}{\partial x_j}=\frac{\partial
f_0}{\partial x_j}\ast f$.

Suppose that $D^\alpha(f_0\ast f)=(D^\alpha f_0)\ast f$. Then
$$D^{\alpha+e_j}(f_0\ast f)=\frac{\partial D^\alpha(f_0\ast f)}{\partial
x_j}=\frac{\partial((D^\alpha f_0)\ast f)}{\partial
x_j}=\frac{\partial D^\alpha f_0}{\partial x_j}\ast
f=(D^{\alpha+e_j}f_0)\ast f.$$ Inductively, we have
$D^\alpha(f_0\ast f)=(D^\alpha f_0)\ast f$ for all multi-index
$\alpha$.

Let $f\in[E^{(\gamma,U)}]$, $\xi\in E$ and $1\leq j\leq n$. By Th.
4.2, $f_0\ast f\in[E^{(\gamma,V)}]$ for some $V\in\mathcal N(E)$ and
$\frac{\partial(f_0\ast f)}{\partial x_j}(\xi)=(f_0\ast
f)(-\frac{\partial\xi}{\partial
x_j})=f[f_0\ast(-\frac{\partial\xi}{\partial
x_j})]=f(-f_0\ast\frac{\partial\xi}{\partial x_j})=f(-\frac{\partial
f_0\ast\xi}{\partial x_j})=\frac{\partial f}{\partial
x_j}(f_0\ast\xi)=(f_0\ast\frac{\partial f}{\partial x_j})(\xi)$,
i.e., $\frac{\partial(f_0\ast f)}{\partial
x_j}=f_0\ast\frac{\partial f}{\partial x_j}$.

If $D^\alpha(f_0\ast f)=f_0\ast D^\alpha f$, then
$D^{\alpha+e_j}(f_0\ast f)=\frac{\partial D^\alpha(f_0\ast
f)}{\partial x_j}=\frac{\partial f_0\ast D^\alpha f}{\partial
x_j}=f_0\ast\frac{\partial D^\alpha f}{\partial x_j} =f_0\ast
D^{\alpha+e_j}f$ since $D^\alpha f\in[E^{(\gamma,W)}]$ for some
$W\in\mathcal N(E)$ by Th. 2.1. $\square$
\end{Theorem}

Recall that if $f\in E'$ for which $supp\,f$ is bounded in $\mathbb
R^n$, then $f$ must be a convolution multiplier [4, Th. 3.3.4]. Then
we can develop the result of continuity of convolution [4, Th.
3.3.5].

First, we give an improvement of Th. 3.3.5 of [4] as follows.

\begin{Theorem}
If $\{f_k\}\subset E'$ such that
$f_k\stackrel{w\ast}{\longrightarrow}f$, i.e., $f_k(\xi)\rightarrow
f(\xi)$ at each $\xi\in E$ ($f\in E'$ by Th. 1.8) and there is a
bounded $F\subset\mathbb R^n$ such that $supp\,f_k\subseteq F,\
\forall\,k\in\mathbb N$, then for every $g\in[E^{(\gamma,U)}]$ and
bounded $B\subset E$,
$$\lim_k(f_k\ast g)(\xi)=(f\ast g)(\xi)\mbox{ uniformly for }\xi\in B.$$

Proof. By Th. 1.9, $supp\,f\subseteq F$ and so $f$ is also a
convolution multiplier on $E$. By Th. 1.2, for every bounded
$B\subset E$, $\lim_kf_k(\xi)=f(\xi)$ uniformly for $\xi\in B$. Then
$f_k\ast\xi\rightarrow f\ast\xi,\ \forall\,\xi\in E$ (see the proof
of Th 3.3.5 of [4]).

Define $T:E\rightarrow E$ and $T_k:E\rightarrow E$ by
$T(\xi)=f\ast\xi$ and $T_k(\xi)=f_k\ast\xi,\ \forall\,\xi\in E,\
k\in\mathbb N$. By Lemma 4.1, $T$ and all $T_k$ are continuous and
linear.

Since $\mathscr D$ is an $(LF)$ space and $\mathscr S$ is a locally
convex Fr\'echet space, both $\mathscr D$ and $\mathscr S$ are
barrelled [3, p.136, 222]. Moreover, $T_k(\xi)=f_k\ast\xi\rightarrow
f\ast\xi=T(\xi)$ at each $\xi\in E$, i.e., $\{T_k(\xi):k\in\mathbb
N\}$ is bounded at each $\xi\in E$. By Th. 9.3.4 of [3], both
$\{T_k:k\in\mathbb N\}$ and $\{T_k:k\in\mathbb N\}\bigcup\{T\}$ are
equicontinuous on $E$.

Let $g\in[E^{(\gamma,U)}]$. Pick a $V\in\mathcal N(E)$ for which
$T(V)\subset U$ and $T_k(V)\subset U$, $\forall\, k\in\mathbb N$,
i.e., $f\ast\eta$, $f_k\ast\eta\in U$, $\forall\,\eta\in V,\
k\in\mathbb N$. Then $f\ast g,\ f_k\ast g\in[E^{(\gamma,V)}]$ for
all $k$ (see the proof of Th. 4.2).

Suppose that $g\in E^{(\gamma,U)}$. Then $f\ast g,\ f_k\ast g\in
E^{(\gamma,V)}$ and
$$(f_k\ast g)(\xi)=g(f_k\ast\xi)\rightarrow g(f\ast\xi)=(f\ast g)(\xi),\ \forall\,\xi\in E.$$
By Th. 1.2, for every bounded $B\subset E$, $\lim_k(f_k\ast
g)(\xi)=(f\ast g)(\xi)$ uniformly for $\xi\in B$.

Now let $g=\sum_{\nu=1}^ma_\nu g_\nu$ where $a_\nu\in\mathbb C$ and
$g_\nu\in E^{(\gamma,U)},\ \nu=1,2,\cdots,m$. Then for every bounded
$B\subset E$ we have that
\begin{align*}
&\lim_k(f_k\ast g)(\xi)=\lim_k
g(f_k\ast\xi)=\lim_k(\sum_{\nu=1}^ma_\nu g_\nu)(f_k\ast\xi)=\lim_k\sum_{\nu=1}^ma_\nu g_\nu(f_k\ast\xi)\\
=&\lim_k\sum_{\nu=1}^ma_\nu(f_k\ast
g_\nu)(\xi)=\sum_{\nu=1}^ma_\nu\lim_k(f_k\ast
g_\nu)(\xi)=\sum_{\nu=1}^ma_\nu(f\ast g_\nu)(\xi)=(f\ast g)(\xi)
\end{align*}
uniformly for $\xi\in B$. $\square$
\end{Theorem}

We also give some simple facts before our main result Th. 4.6.

\begin{Theorem}
Let $f_0\in E'$ be a convolution multiplier, $U\in\mathcal N(E)$ and
$\gamma\in C(0)$. There is a $V\in\mathcal N(E)$ for which
$f_0\ast\cdot:[E^{(\gamma,U)}]\rightarrow[E^{(\gamma,V)}]$ is a
linear operator such that $f_0\ast f\in E^{(\gamma,V)}$ for each
$f\in E^{(\gamma,U)}$. Moreover, if
$f_k\stackrel{w\ast}{\longrightarrow}f$ in $E^{(\gamma,U)}$, i.e.,
$f,f_k\in E^{(\gamma,U)}$ and $f_k(\xi)\rightarrow f(\xi)$ at each
$\xi\in E$, then for every bounded $B\subset E$, $\lim_k(f_0\ast
f_k)(\xi)=(f_0\ast f)(\xi)$ uniformly for $\xi\in B$.

Proof. By Th. 4.2, there is a $V\in\mathcal N(E)$ such that $f_0\ast
f\in E^{(\gamma,V)}$ for $f\in E^{(\gamma, U)}$ and so $f_0\ast
f\in[E^{(\gamma,V)}]$ whenever $f\in[E^{(\gamma,U)}]$.

Let $f,g\in[E^{(\gamma,U)}]$ and $t\in\mathbb C$. Then
$[f_0\ast(f+tg)](\xi)=(f+tg)(f_0\ast\xi)=f(f_0\ast\xi)+tg(f_0\ast\xi)=[(f_0\ast
f)+t(f_0\ast g)](\xi),\ \forall\,\xi\in E$, i.e.,
$f_0\ast(f+tg)=f_0\ast f+t(f_0\ast g)$ and so $f_0\ast\cdot$ is a
linear operator.

Let $f_k\stackrel{w\ast}{\longrightarrow}f$ in $E^{(\gamma,U)}$.
Then $(f_0\ast f_k)(\xi)=f_k(f_0\ast\xi)\rightarrow
f(f_0\ast\xi)=(f_0\ast f)(\xi),\ \forall\,\xi\in E$. By Th. 1.2, for
every bounded $B\subset E$, $\lim_k(f_0\ast f_k)(\xi)=(f_0\ast
f)(\xi)$ uniformly for $\xi\in B$. $\square$
\end{Theorem}

\begin{Corollary}
Suppose that $f_k\stackrel{w\ast}{\longrightarrow}f$ in $E'$ where
each $supp\,f_k\subset\{x\in\mathbb R^n:|x|\leq a\}$ for some $a>0$
and $g_k\stackrel{w\ast}{\longrightarrow}g$ in $E^{(\gamma,U)}$. If
$\xi_k\rightarrow\xi$ in $E$, then
$$(f_k\ast h)(\xi_k)\rightarrow(f\ast h)(\xi),\ \forall\,h\in[E^{(\gamma,U)}],$$
$$(f_0\ast g_k)(\xi_k)\rightarrow(f_0\ast g)(\xi),\ \forall\,\mbox{convolution multiplier }f_0\in E',$$

Proof. Since $\xi_m\rightarrow\xi$, $\{\xi_m:m\in\mathbb N\}$ is
bounded in $E$.

Let $h\in[E^{(\gamma,U)}]$. By Th. 4.4, $\lim_k(f_k\ast
h)(\xi_m)=(f\ast h)(\xi_m)$ uniformly for $m\in\mathbb N$. But
$\lim_m(f_k\ast h)(\xi_m)=(f_k\ast h)(\xi)$ for each $k\in\mathbb N$
and so $\lim_k(f_k\ast
h)(\xi_k)=\lim_{k,m\rightarrow+\infty}(f_k\ast
h)(\xi_m)=\lim_k\lim_m(f_k\ast h)(\xi_m)=\lim_k(f_k\ast
h)(\xi)=(f\ast h)(\xi)$.

Similarly, it follows from Th. 4.5 that for every convolution
multiplier $f_0\in E'$ we have $(f_0\ast
g_k)(\xi_k)\rightarrow(f_0\ast g)(\xi)$. $\square$
\end{Corollary}

\begin{Corollary}
If $f_0\in E'$ is a convolution multiplier and
$P(D)=\sum_{|\alpha|\leq p}a_\alpha D^\alpha$, then
$$P(D)(f_0\ast f)=f_0\ast[P(D)f],\ \forall\,f\in[E^{(\gamma,U)}].$$

Proof. Let $f\in[E^{(\gamma,U)}]$. By Th. 2.1, there is a
$V\in\mathcal N(E)$ such that $D^\alpha f\in[E^{(\gamma,V)}],\
\forall\,|\alpha|\leq p$. By Th. 4.2, there is a $W\in\mathcal N(E)$
such that $f_0\ast D^\alpha f\in[E^{(\gamma,W)}],\
\forall\,|\alpha|\leq p$. By Th. 4.3, Lemma 4.2 and Th. 4.5,
\begin{align*}
P(D)(f_0\ast f)&=\sum_{|\alpha|\leq p}a_\alpha D^\alpha(f_0\ast
f)=\sum_{|\alpha|\leq p}a_\alpha(f_0\ast D^\alpha
f)=\sum_{|\alpha|\leq p}f_0\ast(a_\alpha D^\alpha f)\\
&=f_0\ast(\sum_{|\alpha|\leq p}a_\alpha D^\alpha f)=f_0\ast[P(D)f].\
\square
\end{align*}
\end{Corollary}

We now have a strong continuity result for convolution as follows.

\begin{Theorem}
Let $\{f_k\}\subset E'$ be a sequence of usual distributions such
that $f_k\stackrel{w\ast}{\longrightarrow}f$, i.e.,
$f_k(\xi)\rightarrow f(\xi)$ at each $\xi\in E$ ($f\in E'$ by Th.
1.8) and there is a bounded $F\subset\mathbb R^n$ such that
$supp\,f_k\subseteq F,\ \forall\,k\in\mathbb N$. If
$g_k\stackrel{w\ast}{\longrightarrow}g$ in $E^{(\gamma,U)}$, i.e.,
$g,g_k\in E^{(\gamma,U)}$ for all $k$ and $g_k(\xi)\rightarrow
g(\xi)$ at each $\xi\in E$, then for every bounded $B\subset E$,
$\lim_{k,m\rightarrow+\infty}(f_k\ast g_m)(\xi)=(f\ast g)(\xi)$
uniformly for $\xi\in B$ and, in particular, $\lim_k(f_k\ast
g_k)(\xi)=(f\ast g)(\xi)$ uniformly for $\xi\in B$, and $(f_k\ast
g_k)(\xi_k)\rightarrow(f\ast g)(\xi)$ whenever $\xi_k\rightarrow\xi$
in $E$.

Proof. As in the proof of Th. 4.4, it follows from
$f_k\stackrel{w\ast}{\longrightarrow}f$ in $E'$ and $g,g_k\in
E^{(\gamma,U)}$ that there is a $V\in\mathcal N(E)$ such that $f\ast
g,\ f_m\ast g_k\in E^{(\gamma,V)}$ for all $k,m\in\mathbb N$.

Let $\xi\in E$. As was noticed in the proof of Th. 4.4,
$f_m\ast\xi\rightarrow f\ast\xi$ in $E$ and so $\lim_m(f_m\ast
g_k)(\xi)=\lim_m g_k(f_m\ast\xi)=g_k(f\ast\xi),\
\forall\,k\in\mathbb N$. But $\{f_m\ast\xi\}$ is bounded in $E$ and,
by Th. 1.2, $\lim_k(f_m\ast g_k)(\xi)=\lim_k
g_k(f_m\ast\xi)=g(f_m\ast\xi)=(f_m\ast g)(\xi)$ uniformly for
$m\in\mathbb N$. Then $\lim_{k,m\rightarrow+\infty}(f_m\ast
g_k)(\xi)=\lim_m\lim_k(f_m\ast g_k)(\xi)=\lim_m(f_m\ast
g)(\xi)=\lim_m g(f_m\ast\xi)=g(f\ast\xi)=(f\ast g)(\xi),\
\forall\,\xi\in E$.

Let $B$ be a bounded subset of $E$. If
$\lim_{k,m\rightarrow+\infty}(f_m\ast g_k)(\xi)=(f\ast g)(\xi)$ is
not uniformly for $\xi\in B$, then there exist $\varepsilon>0$,
$\{\xi_\nu\}\subset B$ and integer sequences $k_1<k_2<\cdots$ and
$m_1<m_2<\cdots$ such that
$$|(f_{m_\nu}\ast g_{k_\nu})(\xi_\nu)-(f\ast g)(\xi_\nu)|\geq\varepsilon,\ \nu=1,2,3,\cdots. \leqno(\ast)$$

Since $f\ast g,\ f_{m_\nu}\ast g_{k_\nu}\in E^{(\gamma,V)}$ for all
$\nu\in\mathbb N$ and $$\lim_\nu(f_{m_\nu}\ast
g_{k_\nu})(\xi)=\lim_{k,m\rightarrow+\infty}(f_m\ast
g_k)(\xi)=(f\ast g)(\xi),\ \forall\,\xi\in E,$$ it follows from Th.
1.2 or Th. 1.7 that $\lim_\nu(f_{m_\nu}\ast g_{k_\nu})(\xi)=(f\ast
g)(\xi)$ uniformly for $\xi\in B$ and so there is a $\nu_0\in\mathbb
N$ such that
$$|(f_{m_\nu}\ast g_{k_\nu})(\xi_\nu)-(f\ast g)(\xi_\nu)|<\varepsilon,\ \forall\,\nu>\nu_0.$$
This contradicts ($\ast$) and so
$\lim_{k,m\rightarrow+\infty}(f_m\ast g_k)(\xi)=(f\ast g)(\xi)$
uniformly for $\xi\in B$. $\square$
\end{Theorem}

\end{document}